\documentclass{amsart}
\usepackage{graphicx} 
\usepackage{graphicx} 
\usepackage{amsfonts,amsmath,amssymb}
\usepackage{fullpage}
\usepackage{amsthm}
\usepackage{tikz}
\usetikzlibrary{calc}
\usepackage{mathtools}

\PassOptionsToPackage{obeyspaces}{url}
\usepackage[colorlinks=true,citecolor=blue,urlcolor=blue,linkcolor=blue,bookmarksopen=true]{hyperref}

\theoremstyle{plain}
\newtheorem{theorem}{Theorem}[section]
\newtheorem{lemma}[theorem]{Lemma}
\newtheorem{corollary}[theorem]{Corollary}

\newtheorem{conjecture}[theorem]{Conjecture}
\newtheorem{proposition}[theorem]{Proposition}
\newtheorem{observation}[theorem]{Observation}

\theoremstyle{definition}

\newtheorem{definition}[theorem]{Definition}

\def\P{\mathcal{P}}

\newtheorem{rem}{Remark}

\newcommand{\F}{{\mathcal F}}

\newcommand\cC{{\mathcal C}}
\newcommand\cU{{\mathcal U}}
\newcommand\cD{{\mathcal D}}
\newcommand\bC{{\mathbf C}}
\newcommand\cF{{\mathcal F}}
\newcommand\cE{{\mathcal E}}

\newcommand\cG{{\mathcal G}}
\newcommand\cS{{\mathcal S}}
\newcommand\eps{{\varepsilon}}

\title{On some extremal and probabilistic questions for tree posets}
\author{Bal\'azs Patk\'os and Andrew Treglown}

\thanks{BP: Alfr\'ed R\'enyi Institute of Mathematics, Budapest. Email: \texttt{patkos@renyi.hu}. Research supported by NKFIH grant \\ \indent  FK 132060. \\ \indent
AT: University of Birmingham, United Kingdom. Email: \texttt{a.c.treglown@bham.ac.uk}. Research supported by EPSRC \\ \indent grant EP/V002279/1.
}

\date{}

\date{}
\begin{document}

\maketitle

\begin{abstract}
    Given two posets $P,Q$ we say that $Q$ is $P$-free if $Q$ does not contain a copy of $P$. The size of the largest $P$-free family in $2^{[n]}$, denoted by $La(n,P)$, has been extensively studied since the 1980s. We consider several related problems. For posets $P$ whose Hasse diagrams are trees and have radius at most $2$, we prove that there are 
    $2^{(1+o(1))La(n,P)}$ $P$-free families in $2^{[n]}$, thereby confirming a conjecture of Gerbner,  Nagy, Patk{\'o}s and  Vizer [Electronic Journal of Combinatorics, 2021] in this case. 
    For such $P$ we also resolve the random version of the $P$-free problem, thus generalising the random version of Sperner's theorem due to Balogh, Mycroft and Treglown [Journal of Combinatorial Theory Series A, 2014], and
    Collares Neto and  Morris [Random Structures and Algorithms, 2016]. Additionally, we make a general conjecture that, roughly speaking, asserts that subfamilies of $2^{[n]}$ of size sufficiently above $La(n,P)$ robustly contain $P$, for any poset $P$ whose Hasse diagram is a tree.
\end{abstract}

\section{Introduction}
In this paper we will consider several related extremal and probabilistic questions for posets. All posets  we consider will be  viewed as finite collections of finite subsets of $\mathbb N$ equipped with the containment relation. Indeed,  we will usually be working inside the power set $2^{[n]}$ for some $n \in \mathbb N$.

Let $P,Q$ be posets. A \emph{poset homomorphism} from $P$ to $Q$ is a function $\phi: \ P \rightarrow Q$ such that for every $A,B\in P$, if $A \subseteq B$ then $\phi(A)\subseteq \phi(B)$.
 We say that $P$ is a \emph{subposet} of $Q$ if there is an injective poset homomorphism from $P$ to $Q$; otherwise, $Q$ is said to be \emph{$P$-free}.
The size of the largest $P$-free family in $2^{[n]}$ is denoted by $La(n,P)$. Further we say $P$ is an \emph{induced subposet} of $Q$ if there is an injective poset homomorphism $\phi$ from $P$ to $Q$ such that for every $A,B \in P$, $\phi(A)\subseteq \phi(B)$ if and only if $A\subseteq B$; otherwise, $Q$ is said to be \emph{induced $P$-free}.

The systematic study of $La(n,P)$ was initiated by Katona and Tarj\'an in 1983~\cite{kat}, though Sperner's classical  theorem~\cite{sperner} was the first result in the area. The latter 
asserts that the largest $C_2$-free family\footnote{We write $C_t$ to denote the chain on $t$ elements. Thus, a $C_2$-free family is simply an antichain.} 
in $2^{[n]}$ has size $\binom{n}{\lfloor n/2 \rfloor}$. 
Whilst the asymptotic value of $La(n,P)$ has been determined for a range of posets $P$, the general problem remains open; see~\cite{griggs} or Chapter 7 of~\cite{gerbner2} for a survey on the topic and~\cite{bukh, griggs2} for a conjecture on the asymptotic value of $La(n,P)$ for all posets $P$. 

The \emph{Hasse  diagram} of a poset $P$ is the directed graph with vertex set $P$ and for $A,B \in P$, there is a directed edge from $A$ to $B$ in the Hasse diagram if 
$A \subset B$ and there does not exist a $C \in P$ such that $A \subset C \subset 
B$. 
The \emph{undirected Hasse diagram} is the 
(undirected) graph obtained from the Hasse diagram by removing all orientations on 
the edges. A \emph{tree poset} is a poset whose undirected Hasse diagram is a tree. A tree poset $P$ is  \emph{upward (resp. downward) monotone} if  there exists an element $B \in  P$ with $B \subseteq A$  (resp. $A \subseteq B$) for any $A\in P$. A
tree poset is called \emph{monotone} if it is either upward or downward monotone.

The following result of Bukh~\cite{bukh} determines $La(n,P)$ asymptotically for all tree posets $P$. Given a poset $P$, we write $h(P)$ for the \emph{height of $P$}, i.e., the  length $t$  of  the  longest  chain $C_t$  in $P$.
\begin{theorem}[Bukh~\cite{bukh}]\label{thm:bukh}
    If $P$ is a tree poset then
    $$La(n,P)=(h(P)-1)\binom{n}{\lfloor n/2 \rfloor}(1+O(1/n)).$$
\end{theorem}
A phenomenon often exhibited in combinatorial problems is the property that once one is sufficiently above a `threshold' for guaranteeing the existence of a structure, one can actually `robustly' find this structure. This behaviour is often articulated in terms of \emph{supersaturation}. For example, in recent years there has been significant attention on an old conjecture of Kleitman on the minimum number of chains $C_t$ of length $t$ in a subfamily of $2^{[n]}$ of a given size. Indeed, after earlier progress~\cite{bw, das, dove}, this conjecture was proven by Samotij~\cite{sam}.

The first question we consider takes a different viewpoint on robustness. For this we need the following definition.

\begin{definition}
    Let $P$ be a tree poset and $x\in P$, and let $d\ge 2$ be a positive integer. The \textit{$d$-blow-up $P(x,d)$ rooted at $x$} is the tree poset whose Hasse diagram is defined as follows: for each $u \in P$, if $u$
    is at distance $\rho$ from $x$ in  the undirected Hasse diagram of $P$, then $u$ is replaced with $d^{\rho}$ elements $u^1,u^2,\dots, u^{d^{\rho}}$; furthermore, if $uv$ is an edge of the Hasse diagram of $P$ with $v$ being at distance $\rho-1$ to $x$ in $P$, then the $u^i$s are partitioned into $d^{\rho-1}$ pairwise disjoint sets $U^1,U^2,\dots, U^{d^{\rho-1}}$ each of size $d$, and for every $j \in [d^{\rho-1}]$,
    $v^j$ is joined to all members of $U^j$. The orientation of all such edges is the same as that of $uv$.
\end{definition}

\begin{figure}[ht]\label{pfig}
    \begin{center}
        \pgfdeclarelayer{bg}
        \pgfsetlayers{bg,main}
        \tikzset{vertex/.style={circle, draw=white, thin, fill=black, minimum size = 4pt},
        edge/.style={black, ultra thick}}
        \def\shft{13pt}
        \begin{tikzpicture}
	        \begin{scope}
		        \node [vertex] (p1) at (-3/4,-3/4) {};
    		    \node [vertex] (p2) at (3/4,-3/4) {};
	    	    \node [vertex] (q1) at (-3/4,3/4) {};
		        \node [vertex] (q2) at (3/4,3/4) {};
            \node [vertex] (r) at (3/4,9/4) {};
        		\node [yshift=-\shft] at (p1)   {};
        		\node [yshift=-\shft] at (p2) {};
	        	\node [yshift=\shft] at (q1) {$x$};
	    	    \node [yshift=\shft] at (q2) {};
        		\begin{pgfonlayer}{bg}
	    	        \draw [edge] (p1) -- (q1);
		            \draw [edge] (p2) -- (q1);
		            \draw [edge] (p2) -- (q2);
		            \draw [edge] (q2) -- (r);
          \end{pgfonlayer}
    	    \end{scope}
	
	        \begin{scope}[xshift=80pt]
		        \node [vertex] (p11) at (-3/4,-3/4) {};
    		    \node [vertex] (p12) at (1/4,-3/4) {};
		        \node [vertex] (q1) at (3/4,3/4) {};
          \node [yshift=\shft] at (q1) {$x$};
		        \node [yshift=-\shft] at (p1)   {};
       		\node [yshift=-\shft] at (p2) {};
         \node [vertex] (p21) at (9/4,-3/4) {};
         \node [vertex] (p22) at (17/4,-3/4) {};
         \node [vertex] (q21) at (6/4,3/4) {};
         \node [vertex] (q22) at (12/4,3/4) {};
         \node [vertex] (q23) at (14/4,3/4) {};
         \node [vertex] (q24) at (20/4,3/4) {};
          \node [vertex] (r1) at (5/4,9/4) {};
         \node [vertex] (r2) at (7/4,9/4) {};
         \node [vertex] (r3) at (10/4,9/4) {};
         \node [vertex] (r4) at (12/4,9/4) {};
         \node [vertex] (r5) at (14/4,9/4) {};
         \node [vertex] (r6) at (16/4,9/4) {};
         \node [vertex] (r7) at (20/4,9/4) {};
         \node [vertex] (r8) at (22/4,9/4) {};
        		\begin{pgfonlayer}{bg}
		            \draw [edge] (p11) -- (q1);
		            \draw [edge] (p12) -- (q1);
              \draw [edge] (p21) -- (q1);
		            \draw [edge] (p22) -- (q1);
              \draw [edge] (p21) -- (q21);
		            \draw [edge] (p21) -- (q22);
              \draw [edge] (p22) -- (q23);
		            \draw [edge] (p22) -- (q24);
              \draw [edge] (q21) -- (r1);
		            \draw [edge] (q21) -- (r2);
              \draw [edge] (q22) -- (r3);
		            \draw [edge] (q22) -- (r4);
              \draw [edge] (q23) -- (r5);
		            \draw [edge] (q23) -- (r6);
              \draw [edge] (q24) -- (r7);
		            \draw [edge] (q24) -- (r8);
    	    	\end{pgfonlayer}
    	    \end{scope}
		
	    \end{tikzpicture}
        \caption{A poset rooted at $x$ and its 2-blow-up.}
        \label{fig:Nposet1}
    \end{center}
\end{figure}
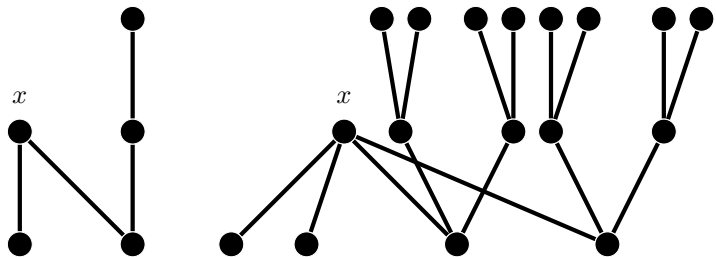

The following conjecture states that once  sufficiently above the threshold in Theorem~\ref{thm:bukh},  one can robustly guarantee a copy of a given tree poset $P$.
\begin{conjecture}\label{conj1}
    Let $P$ be any tree poset of height $h$. Given any $\eps>0$ there exists $\delta >0$ such that the following holds for any sufficiently large $n \in \mathbb N$ and any $x \in P$.
   Consider the $\delta n$-blow-up $P({x, \delta n})$ of $P$  rooted at $x$. If $\mathcal F \subseteq 2^{[n]}$ has size at least $ (h-1+\eps ) \binom{n}{\lfloor n/2 \rfloor}$, then
    $\mathcal F$ contains a copy of $P({x, \delta n})$.
\end{conjecture}

\begin{rem}\label{equiv}
    First observe that, by considering the complement family $\overline{\F}:=\{[n]\setminus F: F\in \F\}$, it is clear that if the conjecture holds for some poset $P$, then it also holds for the dual poset $P^d$ obtained from $P$ by reversing all relations.

    Note that at the cost of decreasing $\delta$ 
    by a constant factor, if Conjecture 
    \ref{conj1} holds for a `simple' tree poset $P$, 
    then it holds for some more complicated tree posets 
    $P'$. More precisely, suppose $u,v$ are leaves 
    of the Hasse diagram of $P'$, $w$ is an 
    element of $P'$ lying on the paths from $x$ to 
    both $u$ and $v$,  and for the unique paths 
    $w=v_0,v_1v_2,\dots,v_k=v$ from $w$ to $v$ and 
    $w=u_0,u_1,u_2,\dots,u_\ell=u$ from $w$ to $u$ 
    with $k\le \ell$ we have $v_{i-1}\subseteq 
    v_{i}$ if and only if $u_{i-1}\subseteq u_{i}$ 
    for all $i\in[k]$.  Finally, suppose that $v_i$ has degree 2 in the Hasse diagram of $P$ for all $i\in [k-1]$. Then if we obtain 
    $P$ from $P'$ by removing $v_1,v_2,\dots, 
    v_k$, then $P'(x,\delta n/2)\subset 
    P(x,\delta n)$ for sufficiently large 
    $n$; so it is enough to verify Conjecture \ref{conj1} for $P$.

    Similarly, suppose $P_1, P_2$ are connected components of $P'\setminus \{x\}$ such that  $y_1\in P_1,y_2\in P_2$ are the only elements of $P_1 \cup P_2$ adjacent to $x$ in the Hasse diagram of $P'$, and also $x\subseteq y_1$ if and only if $x\subseteq y_2$. Suppose further $P_2$ contains a copy of $P_1$ with $y_2$ playing the role of $y_1$. If we obtain $P$ from $P'$ by removing $P_1$, 
    then again $P'(x,\delta n/2)\subset P(x,\delta n)$ for sufficiently large $n$; so it is enough to verify Conjecture~\ref{conj1} for $P$.
\end{rem}

Perhaps the most interesting case of Conjecture~\ref{conj1} is when $P$ is a monotone tree. For example, it implies that if one is sufficiently above the threshold for guaranteeing a binary monotone tree of height $h$ (or even a chain of length $h$), then actually one can find a $\delta n$-ary monotone tree of height $h$.

As well as being interesting in its own right,  Conjecture~\ref{conj1} has implications to a couple of other well-studied questions.
The first of these concerns counting $P$-free families in $2^{[n]}$.
Kleitman~\cite{kleitman} proved that the number of $C_2$-free families in $2^{[n]}$ is $2^{(1+o(1))\binom{n}{\lfloor n/2 \rfloor}}$; this result was further refined by Korshunov~\cite{kor}. More generally, a trivial lower bound on the number of $P$-free families in $2^{[n]}$
is $2^{La(n,P)}$. The following conjecture from~\cite{gerbner} states that this lower bound is close to being tight.
\begin{conjecture}[Gerbner, Nagy, Patk\'os and Vizer~\cite{gerbner}]\label{conj:count}
The number of $P$-free families in   $2^{[n]}$ 
is $$2^{(1+o(1))La(n,P)}.$$
\end{conjecture}

Our first result asserts that if Conjecture~\ref{conj1} is true for a given tree poset $P$, and for a specific choice of some $x \in P$, then
Conjecture~\ref{conj:count} is true for this  $P$.

\begin{theorem}\label{thm:count}
    Suppose $P$ is a tree poset of height $h$ for which Conjecture~\ref{conj1} holds for some choice of $x \in P$. Then the number of $P$-free  families in $2^{[n]}$ is
    $$2^{(h-1+o(1))\binom{n}{\lfloor n/2 \rfloor}}.$$
\end{theorem}
In Section~\ref{sec:count} we prove Theorem~\ref{thm:count} via a variation on the standard \emph{graph container algorithm}.
The following result gives a class of tree posets $P$ for which Conjecture~\ref{conj1} holds for some choices of $x \in P$.

\begin{figure}[ht]\label{posfig}
    \begin{center}
        \pgfdeclarelayer{bg}
        \pgfsetlayers{bg,main}
        \tikzset{vertex/.style={circle, draw=white, thin, fill=black, minimum size = 4pt},
        edge/.style={black, ultra thick}}
        \def\shft{13pt}
        \begin{tikzpicture}
	        \begin{scope}
		        \node [vertex] (p1) at (-3/4,-3/4) {};
    		    \node [vertex] (p2) at (3/4,-3/4) {};
	    	    \node [vertex] (q1) at (-3/4,3/4) {};
		        \node [vertex] (q2) at (3/4,3/4) {};
        		\begin{pgfonlayer}{bg}
	    	        \draw [edge] (p1) -- (q1);
		            \draw [edge] (p2) -- (q1);
		            \draw [edge] (p2) -- (q2);
          \end{pgfonlayer}
    	    \end{scope}
	\begin{scope}[xshift=125pt]
	             \node [vertex] (a) at (0,0) {};
                    \node [vertex] (b1) at (-7/4,5/4) {};
                    \node [vertex] (b2) at (-2/4,5/4) {};
                    \node [vertex] (b3) at (2/4,5/4) {};
                    \node [vertex] (b4) at (7/4,5/4) {};
                    \node [vertex] (c4) at (-4/4,10/4) {};
                    \node [vertex] (c5) at (0,10/4) {};
                    \node [vertex] (c1) at (-11/4,10/4) {};
                    \node [vertex] (c2) at (-9/4,10/4) {};
                    \node [vertex] (c3) at (-7/4,10/4) {};
                    \node [vertex] (c6) at (7/4,10/4) {};

                \begin{pgfonlayer}{bg}
	    	        \draw [edge] (b2) -- (a) -- (b1);
    	    	\draw [edge] (b3) -- (a) -- (b4);
	           \draw [edge] (c4) -- (b2) -- (c5); 
            \draw [edge] (b4) -- (c6);
            \draw [edge] (b1) -- (c1);
            \draw [edge] (b1) -- (c2);
            \draw [edge] (b1) -- (c3);
      \end{pgfonlayer}
              \end{scope}
              
	       \begin{scope}[xshift=275pt]
	             \node [vertex] (a) at (-5/4,0) {};
    		    \node [vertex] (b1) at (-5/4,-5/4) {};
            \node [vertex] (b4) at (0,-5/4) {};
            \node [vertex] (b3) at (-7/4,0) {};
	    	    \node [vertex] (b2) at (0,0) {};
		        \node [vertex] (c) at (5/4,0) {};
                \node [vertex] (d) at (0,5/4) {};
                \node [vertex] (e) at (0,10/4) {};
                \node[vertex] (f) at (-5/4,-10/4) {};
                \node[vertex] (f2) at (-2/4,-10/4) {};
                \node[vertex] (f5) at (2/4,-10/4) {};
                 \node[vertex] (f3) at (-7/4,-10/4) {};
                  \node[vertex] (f4) at (-9/4,-10/4) {};

        		\begin{pgfonlayer}{bg}
	    	        \draw [edge] (a) -- (b1) -- (b2) -- (d) -- (c);
    	    	\draw [edge] (f) -- (b1);
	    	\draw [edge] (e) -- (d);        \draw [edge] (e) -- (d) -- (b2); 
            \draw [edge] (b1) -- (b3);      
            \draw [edge] (f3) -- (b1) -- (f4);
            \draw [edge] (d) -- (b4) -- (f2);
            \draw [edge] (f5) -- (b4);
      \end{pgfonlayer}
    	    \end{scope}
		
	    \end{tikzpicture}
        \caption{Some posets of radius 2: the $N$ poset, a monotone tree, and a more complicated one.}
        \label{fig:Nposet}
    \end{center}
\end{figure}
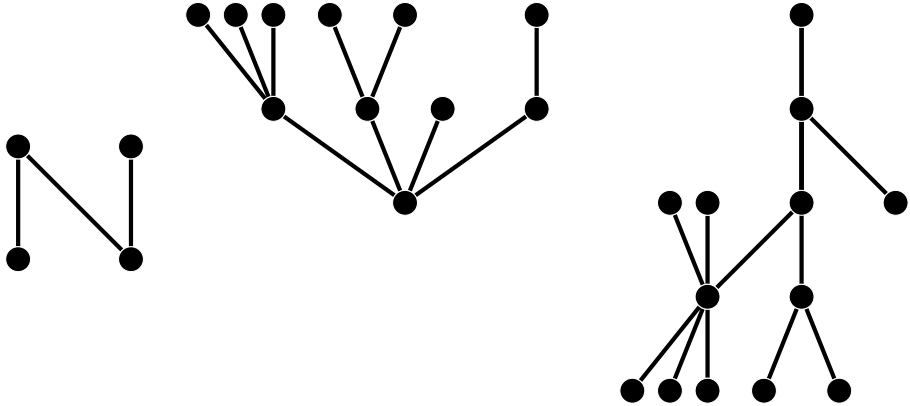

\begin{theorem}\label{thmtrue}
     Conjecture~\ref{conj1} holds for 
     all tree posets $P$ and any $x\in P$ for which all elements of $P$ are within distance at most 2 of $x$ in the (undirected) Hasse diagram of $P$. 
\end{theorem}
In Section~\ref{sec:further} we discuss why it seems challenging to extend Theorem~\ref{thmtrue} to  arbitrary tree posets.
 The next corollary follows immediately from Theorems~\ref{thm:count} and~\ref{thmtrue}.
Here we say a tree poset $P$ has \emph{radius at most $t$} if there is some $x\in P$ so that every vertex in the (undirected) Hasse diagram is at distance at most $t$ from $x$.
\begin{corollary}\label{cor1}
    Given any tree poset $P$ of height $h\leq 5$ and radius at most $2$, the number of $P$-free  families in $2^{[n]}$ is
    $$2^{(h-1+o(1))\binom{n}{\lfloor n/2 \rfloor}}.$$
\end{corollary}
Notice that every tree poset $P$ of radius at most $2$ has height at most $5$; we only
include the $h\leq 5$ condition in Corollary~\ref{cor1} (and in the statement of other results below) to make this restriction on the height explicit.
Note that Corollary~\ref{cor1} resolves Conjecture~\ref{conj:count} for various natural posets $P$ including all monotone trees of height at most $3$ and the poset $N$ (see Figure~\ref{fig:Nposet}).

We also prove the following special version of Theorem~\ref{thmtrue} that ensures a larger blow-up of $P$. 
\begin{lemma}\label{lem:special}
 Let $P$ be any tree poset of height $h \leq 5$ and radius at most $2$. Let $x \in P$ such that every element of $P$ is of distance at most two in the (undirected) Hasse diagram of $P$.
Given any $\eps>0$, the following holds for any sufficiently large $n \in \mathbb N$.
Consider the $ n^{1.9}$-blow-up $P({x,  n^{1.9}})$ of $P$ rooted at $x$. If $\mathcal F \subseteq 2^{[n]}$ is such that $|\mathcal F| \geq 4(h-1+\eps ) \binom{n}{\lfloor n/2 \rfloor}$, then
    $\mathcal F$ contains a copy of $P({x,  n^{1.9}})$.  
\end{lemma}
This lemma is a crucial tool in the proof of Theorem~\ref{thm:randomgeneral}.
Note that we have not tried to optimise the lower bound on $|\cF|$ in Lemma~\ref{lem:special} as the result suffices for our purposes.


Our other application of Conjecture~\ref{conj1} is to the random version of the $P$-free problem. Let $\P (n,p)$ be obtained from $2^{[n]}$ by selecting each element of  $2^{[n]}$ independently at random with probability $p$. The model $\P(n,p)$ was first investigated by R\'enyi~\cite{renyi} who determined the probability threshold for the property that $\P(n,p)$ is not itself an antichain, thereby answering a question of Erd\H{o}s. There has also been interest in the problem of determining when the largest $P$-free family of $\P (n,p)$ has size $p \cdot (1+o(1)) La(n,P)$ with high probability (w.h.p.).\footnote{Here, by \emph{with high probability} we mean with probability tending to $1$ as $n$ tends to infinity.} 
For example, the next result resolves this problem when $P=C_2$.

\begin{theorem}[Balogh, Mycroft and Treglown~\cite{bmt} and Collares Neto and  Morris~\cite{cm}] \label{randomsperner}
For any $\eps > 0$ there exists a constant $C$ such that if $p > C/n$ then w.h.p. the largest $C_2$-free family (i.e., antichain) in $\P(n,p)$ has size  $(1\pm\eps)p \binom{n}{\lfloor n/2 \rfloor}$.
\end{theorem}
Osthus~\cite{osthus} had earlier proven Theorem~\ref{randomsperner} in the case when $pn/\log n \rightarrow \infty$. Moreover, Osthus showed that, for a fixed $c>0$, if $p=c/n$ then with high probability the largest antichain in $\P(n,p)$ has size at least $(1+o(1))(1+e^{-c/2}) p \binom{n}{\lfloor n/2 \rfloor}$. So the bound on $p$ in Theorem~\ref{randomsperner} is best-possible up to the constant $C$. 

The results in~\cite{bmt, cm} actually yield an analogous result for $C_t$-free families (for any fixed $t \in \mathbb N$). Hogenson~\cite{hog} also adapted the proof of Theorem~\ref{randomsperner} from~\cite{bmt} to obtain an analogous result for the so-called $t$-fork. In~\cite[Conjecture~7]{gerbner} a general conjecture was made on the threshold for the $P$-free problem in $\P(n,p)$ and the corresponding problem for induced $P$ was also considered. Our next result resolves the random version of the $P$-free problem for all tree posets of radius at most $2$.

\begin{theorem}\label{thm:randomgeneral}
Let $P$ be any tree poset of height $h \leq 5$ and radius at most $2$.
Given any $\eps>0$, there exists $C=C(\eps, P)>0$ such that the following holds. If $p>C/n$ then w.h.p. the largest $P$-free family in $\mathcal P(n,p) $ has size 
$(h- 1\pm \eps )p\binom{n}{\lfloor n/2 \rfloor}.$ 
\end{theorem}
Notice that Theorem~\ref{thm:randomgeneral} implies the random version of Sperner's theorem (Theorem~\ref{randomsperner}).
Furthermore, the lower bound on $p$ in Theorem~\ref{thm:randomgeneral} is essentially tight. Indeed, for such $P$, \cite[Corollary~6]{gerbner} implies that if $p=o(1/n)$ then w.h.p. the largest $P$-free family in $\mathcal P(n,p) $ has size at least 
$(h-o(1) )p\binom{n}{\lfloor n/2 \rfloor}.$
The proof of Theorem~\ref{thm:randomgeneral} applies Theorem~\ref{thmtrue}, Lemma~\ref{lem:special} and the graph container method. As noted by one of the referees of the paper, it may be possible to prove Theorem~\ref{thm:randomgeneral} via the \emph{hypergraph} container method. In fact, the analogous result for $C_t$-free families
was proven  using the {hypergraph} container method in~\cite{cm}.
Our proofs of Theorems~\ref{thm:count} 
and~\ref{thm:randomgeneral} show that variants of the  graph container algorithm have the potential to attack problems that cannot be naturally stated in the language of independent sets in graphs. This philosophy is perhaps  the most important message for the reader to take from our work.

\smallskip 

The paper is organised as follows. In Section~\ref{sec:pre} we introduce some useful results. The proofs of Theorems~\ref{thm:count}, \ref{thmtrue} and~\ref{thm:randomgeneral} and Lemma~\ref{lem:special} are presented in Section~\ref{sec:proofs}.
In Section~\ref{sec:further} we provide some further discussion and make a few simple observations.

\smallskip

{\noindent \bf Notation.} Throughout this paper we omit floors and ceilings whenever this
does not affect the argument.
Given a set $S$ and $t \in \mathbb N$ we write $\binom{S}{\leq t}$ to denote the set of all subsets of $S$ of size at most $t$, and let $\binom{|S|}{\leq t}:=\left |\binom{S}{\leq t}\right |$.

Given any poset $P$, when considering $P$ as a Hasse diagram we will often refer to the elements of $P$ as vertices.
Consider a tree poset $P$ with root $x \in P$. 
Suppose $u \in P$ is of distance $\rho$ from $x$ in the undirected Hasse diagram of $P$, and let $v \in P$ be a neighbour of $u$ of   distance $\rho+1$ from $x$ 
in the undirected Hasse diagram. Then we call $u$ the \emph{parent of $v$} and $v$ a \emph{child of $u$}. Note that every vertex other than $x$ has precisely one parent, though every vertex can have multiple children.

\section{Preliminaries}\label{sec:pre}

In this section, we introduce a few tools  that will be used in the proofs of some of our main results.

\begin{definition}
    The \textit{Lubell-mass} $\lambda_n(\cF)$ of a family $\cF\subseteq 2^{[n]}$ is $\sum_{F\in\cF}\frac{1}{\binom{n}{|F|}}$. This is the expected number of sets in $\cF$ that are contained in a maximum chain $\cC$ in $2^{[n]}$ chosen uniformly at random among all such chains.

    Let $\bC_n$ denote the set of all maximal chains in $2^{[n]}$. Given a family $\cF\subseteq 2^{[n]}$, the\textit{ min-partition} of $\bC_n$ is $\bigcup_{F\in \cF}\bC_F \cup \bC ^*_\emptyset$, where $\bC_F$ contains those $\cC\in \bC_n$ where $F$ is the smallest set in $\cF\cap \cC$, while $\bC_\emptyset ^*$ is the set of those $\cC\in \bC_n$ for which $\cC\cap \cF=\emptyset$.
    For a set $F\in \cF$, we write $\cU_\cF(F):=\{G\setminus F:\, G \in  \cF \, , F\subseteq G \}$ and we drop the subscript $\cF$ if it is clear from context.
\end{definition}

As $\lambda_n(\cF)$ is the expected number of sets of $\cF$ in a randomly selected $\cC\in \bC_n$, and $\lambda_{n-|F|}(\cU(F))$ is the expected number of sets of $\cF$ in a randomly selected chain $\cC\in \bC_F$, we obtain the following proposition. 

\begin{proposition}\label{average}
    Given $K>0$, if $\lambda_{n-|F|}(\cU(F))\le K$ for all $F\in \cF$, then $\lambda_n(\cF)\le K$.
\end{proposition}

\begin{proposition}[Griggs, Li and Lu, Lemma 3.2 in  \cite{griggs3}]\label{size}
   Given $K>0$,  if $\lambda_n(\cF)\le K$, then $|\cF|\le K\binom{n}{\lfloor n/2\rfloor}$.
\end{proposition}

The proof of the next lemma is essentially identical to that of Theorem 3.1 in \cite{griggs4}. Let $\nabla^j_\cF(F):=|\{G\in \cF: F\subseteq G, |G|=|F|+j\}|$.

\begin{lemma}\label{fork}
   Let $\eps>0$ and  $n \in \mathbb N$ be sufficiently large. Let $\cF \subseteq 2^{[n]}$ be a family that contains only sets of size in $[n/2-n^{2/3},n/2+n^{2/3}]$. If for every $F\in \cF$ we have (i) $\nabla^1_\cF(F)\le \varepsilon n$ and (ii) $\nabla^j_\cF(F) \le \varepsilon n^2$ for all $j\ge 2$, then $|\cF|\le (1+15\varepsilon)\binom{n}{\lfloor n/2\rfloor}$.
\end{lemma}

\begin{proof}
    We consider the min-partition of $\bC_n$ and bound $\lambda_{n-|F|}(\cU(F))$. As all sets of $\cF$ have size in $[n/2-n^{2/3},n/2+n^{2/3}]$, we obtain that for all $F \in \mathcal F$,
    $$\lambda_{n-|F|}(\cU(F))\le 1+\frac{\nabla^1_\cF(F)}{n-|F|}+\frac{\nabla^2_\cF(F)}{\binom{n-|F|}{2}}+\sum_{j=3}^{2n^{2/3}} \frac{\nabla^j_\cF(F)}{\binom{n-|F|}{j}}\le 1+14\varepsilon+O\left (\frac{\varepsilon}{n} \right),$$
    where here we used properties (i) and (ii). Propositions~\ref{average} and~\ref{size} complete the proof.
\end{proof}

\begin{lemma}\label{fork+}
    Let $\eps>0$ and  $n \in \mathbb N$ be sufficiently large.
    Let $\cF \subseteq 2^{[n]}$ be a family that contains only sets of size in $[n/2-n^{2/3},n/2+n^{2/3}]$. If for every $F\in \cF$ we have  $\nabla^j_\cF(F) \le \varepsilon n^4$ for all $j\ge 4$, then $|\cF|\le (4+400\varepsilon)\binom{n}{\lfloor n/2\rfloor}$.
\end{lemma}

\begin{proof}
    We again consider the min-partition of $\bC_n$ and bound $\lambda_{n-|F|}(\cU(F))$. For
     all $F \in \mathcal F$,
    $$\lambda_{n-|F|}(\cU(F))\le 4+\frac{\nabla^4_\cF(F)}{\binom{n-|F|}{4}}+\sum_{j=5}^{2n^{2/3}} \frac{\nabla^j_\cF(F)}{\binom{n-|F|}{j}}\le 4+399\varepsilon+O\left (\frac{\varepsilon}{n}\right ).$$
     Propositions~\ref{average} and~\ref{size} complete the proof.
\end{proof}

\section{The proofs of our main results}\label{sec:proofs}
\subsection{Proof of Theorem~\ref{thm:count}}\label{sec:count}
To prove Theorem~\ref{thm:count} we apply the following graph container result.
\begin{lemma}\label{lem:con}
Suppose $P$ is a tree poset of height $h$ for which Conjecture~\ref{conj1} holds for some choice of $x \in P$.  
Given any $\eps>0$, let $\delta >0$ be as in Conjecture~\ref{conj1} and suppose $n \in \mathbb N$ is sufficiently large. 
Then there exists a function $f : \binom{2^{[n]}}{\leq |P| 2^n/\delta n} \to \binom{2^{[n]}}{\leq (h-1+\eps)\binom{n}{\lfloor n/2\rfloor }}$  such that 
for any $P$-free family $\mathcal F$ in $2^{[n]}$ there is a subfamily $\mathcal H \subseteq \mathcal F$ so that $\mathcal H \in \binom{2^{[n]}}{\leq |P| 2^n/\delta n}$ and $\mathcal F \subseteq \mathcal H \cup f(\mathcal H)$.
\end{lemma}
Given any $\mathcal H \in \binom{2^{[n]}}{\leq |P| 2^n/\delta n}$ we refer to the family $\mathcal H \cup f(\mathcal H)$ produced by Lemma~\ref{lem:con} as a \emph{container}. We call $\mathcal H$ the \emph{fingerprint} of the container 
$\mathcal H \cup f(\mathcal H)$.
With Lemma~\ref{lem:con} at hand, Theorem~\ref{thm:count} now follows easily.

\begin{proof}[Proof of Theorem~\ref{thm:count}]
Suppose $P$ is a tree poset of height $h$ for which Conjecture~\ref{conj1} holds for some choice of $x \in P$.  By Theorem~\ref{thm:bukh}, there exists a $P$-free family $\mathcal F$ in $2^{[n]}$ of size $(h-1+o(1))\binom{n}{\lfloor n/2 \rfloor}$.
By considering all subfamilies of $\mathcal F$ we see that there are at least 
$2^{ (h-1+o(1))\binom{n}{\lfloor n/2 \rfloor} }$ $P$-free families in $2^{[n]}$.
To prove the theorem, note that it therefore suffices to prove the following: given any $\eps >0$, if $n\in \mathbb N$ is sufficiently large then there are at most 
$2^{ (h-1+3\eps)\binom{n}{\lfloor n/2 \rfloor} }$ $P$-free families in $2^{[n]}$.

Given any $\eps>0$ let  $\delta >0$ be as in Conjecture~\ref{conj1} and suppose $n \in \mathbb N$ is sufficiently large. Let $\mathcal C$ be the collection of all the containers $\mathcal H \cup f(\mathcal H)$ (for all $\mathcal H \in \binom{2^{[n]}}{\leq |P| 2^n/\delta n}$) obtained when applying Lemma~\ref{lem:con}.
Thus,
$$|\mathcal C| = \binom{2^{n}}{\leq |P| 2^n/\delta n}\leq 2^{\frac{2|P|2^{n}}{\delta n} \log n}\leq 2^{\eps \binom{n}{\lfloor n/2 \rfloor}},$$
where the last inequality follows as $n$ is sufficiently large.
Moreover, by Lemma~\ref{lem:con}, $|C|\leq (h-1+2 \eps )\binom{n}{\lfloor n/2 \rfloor}$ for every $C \in \mathcal C$. Since every $P$-free family in $2^{[n]}$ is a subset of some container $C \in \mathcal C$, we immediately obtain that there are at most
$$|\mathcal C| \cdot 2^{ (h-1+2\eps)\binom{n}{\lfloor n/2 \rfloor}}\leq 2^{ (h-1+3\eps)\binom{n}{\lfloor n/2 \rfloor}}$$
$P$-free families in $2^{[n]}$, as desired.
\end{proof}
It remains to prove Lemma~\ref{lem:con}. The proof uses a modified version of the graph container algorithm. In particular,  there are some specific subtleties we need to take care of.

In the proof of Lemma~\ref{lem:con} we will be constructing our containers through a process of analysing copies of $t$-blow-ups $P(x,t)$ in $2^{[n]}$. For this, it will be helpful to consider the following process of `constructing' a copy of $P$ in $P(x,t)$ with root $x$. Suppose the children of $x$ in $P$ are
$u_1, \dots, u_{s}$. Then for each of these children $u_i$ there is a set $U_i$ of size $t$ in $P(x,t)$ where each $x_i \in U_i$ can `play the role' of $u_i$ in a copy of $P$ rooted at $x$. That is, there is an edge between $x_i $ and $x$ in the Hasse diagram of $P(x,t)$ and it is oriented in the `same' direction as the edge between $u_i $ and $x$ in the Hasse diagram of $P$.

Fix some choice of  $x_i \in U_i$ for each $i \in [s]$; so $x_i$ will play the role of $u_i$ in the copy of $P$ we are constructing. 
Fix  any $i \in [s]$, and let $u_{1,i},\dots, u_{s_i,i}$ denote the children of $u_i$ in $P$. Again, for each of these children $u_{j,i}$ there is a set $U_{j,i}$ of size $t$ in $P(x,t)$ where each $x_{j,i} \in U_{j,i}$ can play the role of $u_{j,i}$ in a copy of $P$ rooted at $x$ that contains $x_1,\dots, x_s$. 
In particular, every vertex $x_{j,i} \in U_{j,i}$ is adjacent to $x_i$ in the 
Hasse diagram of $P(x,t)$ and the edge between $x_{j,i}$ and $x_i$ is oriented in the `same' direction as the edge between $u_{j,i} $ and $u_i$ in the Hasse diagram of $P$.
For each $i \in [s]$ and $j \in [s_i]$, fix some  $x_{j,i} \in U_{j,i}$ to play the role of $u_{j,i}$ in the copy of $P$ we are constructing.
We repeat this process to construct a copy of $P$ rooted at $x$ in $P(x,t)$: at every step we consider the children $u$ of every vertex selected in the previous step;  we have $t$ vertices in $P(x,t)$ that are able to play the role of $u$ in the copy of $P$ we are constructing, so we can continue this process to obtain our desired copy of $P$.

Suppose now  $\mathcal F$ is a $P$-free family in $2^{[n]}$ and consider a copy of 
$P(x,t)$ in $2^{[n]}$ for some $t \in \mathbb N$ and $x \in P$. Run the process described in the previous two paragraphs, however, now at every step insist that any vertex in $P(x,t)$ selected must be in $\mathcal F$. 
Since $\mathcal F$ is a $P$-free family this process cannot finish and produce a copy of $P$. That is, either the root of our copy of 
$P(x,t)$ does not belong to 
$ \mathcal F$, or at some step of the process we must encounter a set $U$ of $t$ vertices that corresponds to some $u \in P$ and where $U \cap \mathcal F = \emptyset$. 
As we will now see, this simple observation is crucial to the success of running our version of  the graph container algorithm.

\begin{proof}[Proof of Lemma~\ref{lem:con}]
Suppose $P$ is  a tree poset of height $h$ for which Conjecture~\ref{conj1} holds for some choice of $x \in P$.  Set $p:=|P|$.
Given any $\eps>0$, let $\delta >0$ be as in Conjecture~\ref{conj1} and suppose $n \in \mathbb N$ is sufficiently large. 

Fix a total order $\mathcal O_{2^{[n]}}$ of the elements of ${2^{[n]}}$.
Let $\mathcal O_{P}$ be a total order of the vertices of $P$ such that the first vertex is $x$; the next vertices are the children of $x$; the next vertices are those vertices of distance two from $x$ in the undirected Hasse diagram of $P$, and so forth. 
We write $y_1,\dots, y_{p}$ for the vertices of $P$ ordered as in $\mathcal O_{P}$; so $y_1=x$.
Further, let $\mathcal P_{\text{blow}}$ denote the set of all copies of $P(x,t)$ in $2^{[n]}$ for all $t \in \mathbb N$. Let $\mathcal O_{\text{blow}}$ be a total order of the elements of $\mathcal P_{\text{blow}}$.

We now run our modified version of the graph container algorithm. The input of the algorithm is a $P$-free family $\mathcal F \subseteq 2^{[n]}$. The algorithm will output a fingerprint $\mathcal H$ and a container $\mathcal H \cup f(\mathcal H)$ where 
$\mathcal H \subseteq \mathcal F \subseteq \mathcal H \cup f(\mathcal H)$. 

Initially we set $\mathcal G^0:= 2^{[n]}$ and $\mathcal H^0:=\emptyset$. We will add vertices from $2^{[n]}$ to $\mathcal H^0$ and remove vertices from $\mathcal G^0$ through the following iterative process, beginning at Step~$1$.

At Step~$i$, let $P_1(x,t) \in \mathcal P_{\text{blow}}$ be a copy of $P(x,t)$ in 
$\mathcal G^{i-1}\subseteq 2^{[n]}$ where we choose $t$ to be as large as possible. If there is more than one copy of $P(x,t)$ in $\mathcal G^{i-1}$ for this choice of $t$ we choose $P_1(x,t)$ to be the copy of $P(x,t)$ appearing earliest in the total order $\mathcal O_{\text{blow}}$. Let $x_1$ be the vertex of $P_1(x,t)$ that plays the role of $x$. 

\begin{itemize}
    \item { Suppose $x_1 \not \in \mathcal F$.} Then define 
    $\mathcal G^{i}:= \mathcal G^{i-1} \setminus \{x_1\}$ and $\mathcal H^i:=\mathcal H^{i-1}$. Proceed to Step~$i+1$.

    \item Suppose $x_1  \in \mathcal F$ and $t \geq \delta n$. We will update 
    $\mathcal G^{i-1}$ and $\mathcal H^{i-1}$ in several \emph{subphases} to obtain
    $\mathcal G^{i}$ and $\mathcal H^{i}$ respectively.
\begin{itemize}
 \item {\it Subphase~$1$.}   Delete $x_1$ from $\mathcal G^{i-1}$ and add it to $\mathcal H^{i-1}$. Proceed to Subphase~$2$.

 \item {\it Subphase~$2$.} 
 Recall that we write $y_1,\dots, y_{p}$ for the vertices of $P$ ordered as in $\mathcal O_{P}$. So $y_1=x$ and $y_2$ is a child of $x$ in $P$. In $P_1(x,t)$ there is a set of $t$ vertices $U_2 \subseteq \mathcal G^{i-1}$ so that each $x _2 \in U_2$ is adjacent to $x_1$ in the Hasse diagram of $P_1(x,t)$ and further the edge between
$x_1$ and $x_2$ is oriented in the same direction as the edge between $y_1(=x)$ and $y_2$ in the Hasse diagram for $P$. 

\smallskip

Look through the vertices  in $U_2\subseteq 2^{[n]}$ one by one, following the total order $\mathcal O_{2^{[n]}}$. If we run through the entire set $U_2$ without finding an element from $\mathcal F$ then we delete all $t$ vertices in $U_2$ from 
$\mathcal G^{i-1}$ and stop Step~$i$.

\smallskip

Otherwise, let $x_2 \in U_2$ denote the very first element from $\mathcal F$ that we discover in $U_2$. Add $x_2$ to $\mathcal H^{i-1}$ and remove from 
$\mathcal G^{i-1}$ both  $x_2$ and all elements of $U_2$ that occur before $x_2$ in the total order
$\mathcal O_{2^{[n]}}$. So here we have deleted between $1$ and $t$ elements of
$\mathcal G^{i-1}$, and of these  only $x_2$ lies in $\mathcal F$. Proceed to Subphase~$3$.

\item {\it Subphase~$j$ (for $3\leq j \leq p$).} 
In the previous subphases we have defined $x_1,\dots, x_{j-1} \in P_1(x,t)$ corresponding to $y_1,\dots, y_{j-1} \in P$.

\smallskip

By definition of $\mathcal O_{P}$, the unique parent $y_k$ of $y_j$ in $P$ must be one of  $y_1,\dots, y_{j-1}$. Furthermore, by definition of $P(x,t)$,
there is a set of $t$ vertices $U_j\subseteq \mathcal G^{i-1}$ in $P_1(x,t)$ corresponding to $y_j$
so that  every  $x _j \in U_j$ is adjacent to $x_k$ in the Hasse diagram of $P_1(x,t)$, and moreover, the edge between
$x_k$ and $x_j$ is oriented in the same direction as the edge between $y_k$ and $y_j$ in the Hasse diagram for $P$. 

\smallskip

Look through the vertices  in $U_j\subseteq 2^{[n]}$ one by one, following the total order $\mathcal O_{2^{[n]}}$. If we run through the entire set $U_j$ without finding an element from $\mathcal F$ then we delete all $t$ vertices in $U_j$ from 
$\mathcal G^{i-1}$ and stop Step~$i$.

\smallskip

Otherwise, let $x_j \in U_j$ denote the very first element from $\mathcal F$ that we discover in $U_j$. Add $x_j$ to $\mathcal H^{i-1}$ and remove from 
$\mathcal G^{i-1}$ both $x_j$ and all elements of $U_j$ that occur before $x_j$ in the total order
$\mathcal O_{2^{[n]}}$. So here we have deleted at most $t$ elements of
$\mathcal G^{i-1}$, of which only $x_j$ lies in $\mathcal F$. Proceed to Subphase~$j+1$.

\end{itemize}

\smallskip

At the end of these subphases relabel $\mathcal G^{i-1}$ as 
$\mathcal G^{i}$ and $\mathcal H^{i-1}$ as $\mathcal H^i$; then proceed to Step~$i+1$. 
Since $\mathcal F$ is $P$-free, for some $j \in [p]$, Subphase~$j$ must consider a set $U_j$ where $U_j \cap \mathcal F =\emptyset$. Thus, 
 in Step~$i$ we have deleted at least $t$ vertices from $\mathcal G^{i-1}$ to obtain $\mathcal G^i$ and have added at most $p-1$ vertices to 
$\mathcal H^{i-1}$ (at most one in each of the first $p-1$ subphases) to obtain $\mathcal H^i$.

    \item Suppose $x_1  \in \mathcal F$ and $t < \delta n$. Then define
     $\mathcal H:=\mathcal H^{i-1} \cup \{x_1\}$ and $f(\mathcal H):= \mathcal G^{i-1} \setminus \{x_1\}$ and terminate the algorithm.
\end{itemize}

Note that in every step of the algorithm we only add elements to $\mathcal H^{i-1}$ that lie in $\mathcal F$. So certainly $\mathcal H \subseteq \mathcal F$. Similarly, by construction, $\mathcal F \subseteq \mathcal H \cup f(\mathcal H)$.
In every step of the algorithm, except the final step, if we add (at most $p-1$) elements to $\mathcal H^{i-1}$ we delete at least $\delta n$ elements from 
$\mathcal G^{i-1} \subseteq 2^{[n]}$.  Therefore,
$$| \mathcal H|\leq \frac{(p-1)2^n}{ \delta n}+1 \leq \frac{p \cdot 2^n}{ \delta n}. $$
Furthermore, by construction of the algorithm, the family $f(\mathcal H)\subseteq 2^{[n]}$ does not contain a copy of the $\delta n$ blow-up $P(x,\delta n)$. Thus, by the assumption of the lemma, we have that 
$$|f(\mathcal H)|\leq (h-1+\eps ) \binom{n}{\lfloor n/2 \rfloor}.$$

Therefore all that remains to check is that the function $f$ is well-defined. That is, suppose $\mathcal F_1$ and $\mathcal F_2$ are $P$-free families in $2^{[n]}$
such that, on input $\mathcal F_i$, our  algorithm outputs 
the container $\mathcal H_i \cup f(\mathcal H_i)$ for each $i \in [2]$;
if $\mathcal H_1=\mathcal H_2$ then we require that $f(\mathcal H_1)=f(\mathcal H_2)$. This follows though immediately from the definition of the algorithm. 
Indeed, suppose one is only presented with the fingerprint that is outputted by the algorithm. Then one can completely identify every action taken during the algorithm. Thus,
two applications of the algorithm yielding the same fingerprint must also yield the same container. 
\end{proof}

\subsection{Proofs of Theorem~\ref{thmtrue} and Lemma \ref{lem:special}}

\begin{figure}[ht]
    \begin{center}
        \pgfdeclarelayer{bg}
        \pgfsetlayers{bg,main}
        \tikzset{vertex/.style={circle, draw=white, thin, fill=black, minimum size = 4pt},
        edge/.style={black, ultra thick}}
        \def\shft{13pt}
        \begin{tikzpicture}
	        \begin{scope}
    		    \node [vertex] (b) at (-3/4,-3/4) {};
	    	    \node [vertex] (c) at (0,3/4) {};
		        \node [vertex] (d) at (3/4,-3/4) {};
        		\node [yshift=-\shft] at (b) {$x$};
            \node [yshift=\shft] at (c) {$c$};
            \node [yshift=-\shft] at (d) {$a$};
	    
        		\begin{pgfonlayer}{bg}
	    	        \draw [edge] (b) -- (c) -- (d);
    	    	\end{pgfonlayer}
    	    \end{scope}
	
	        \begin{scope}[xshift=100pt]
	           \node [vertex] (a) at (-3/2,-3/4) {};
    		    \node [vertex] (b) at (-3/4,3/4) {};
	    	    \node [vertex] (c) at (0,-3/4) {};
		        \node [vertex] (e) at (-3/4,9/4) {};
        
	        	\node [yshift=-\shft] at (c) {$x$};
                    \node [yshift=-\shft] at (a) {$a$};
                    \node [xshift=-\shft] at (b) {$c$};
                    \node [xshift=-\shft] at (e) {$d$};
	    	    
        		\begin{pgfonlayer}{bg}
	    	        \draw [edge] (a) -- (b) -- (c);
              \draw [edge] (b) --  (e);
    	    	\end{pgfonlayer}
    	    \end{scope}
    	    
    	   \begin{scope}[xshift=175pt]
	           \node [vertex] (a) at (-5/4,0) {};
    		    \node [vertex] (b1) at (0,-5/4) {};
	    	    \node [vertex] (b2) at (0,0) {};
		        \node [vertex] (b3) at (0,5/4) {};
		        \node [vertex] (c) at (5/4,0) {};
        	
	        	\node [xshift=-\shft] at (b2) {$x$};
	        	\node [xshift=\shft] at (b1) {$b$};
                    \node [xshift=-\shft] at (b3) {$d$};
                    \node [yshift=\shft] at (a) {$a$};
                    \node [yshift=-\shft] at (c) {$e$};
        		\begin{pgfonlayer}{bg}
	    	        \draw [edge] (a) -- (b1) -- (b2) -- (b3) -- (c);
    	    	\end{pgfonlayer}
    	    \end{scope}
         \begin{scope}[xshift=275pt]
	           \node [vertex] (a) at (-5/4,0) {};
    		    \node [vertex] (b1) at (0,-5/4) {};
	    	    \node [vertex] (b2) at (0,0) {};
		        \node [vertex] (c) at (5/4,0) {};
                \node [vertex] (d) at (0,5/4) {};
                \node [vertex] (e) at (0,10/4) {};
        	
	        	\node [xshift=-\shft] at (b2) {$x$};
                    \node [xshift=\shft] at (b1) {$b$};
                    \node [xshift=-\shft] at (d) {$d$};
                    \node [xshift=-\shft] at (e) {$f$};
                    \node [yshift=\shft] at (a) {$a$};
                    \node [yshift=-\shft] at (c) {$e$};
	        	
        		\begin{pgfonlayer}{bg}
	    	        \draw [edge] (a) -- (b1) -- (b2) -- (d) -- (c);
    	    	
	    	        \draw [edge] (e) -- (d) -- (b2); \end{pgfonlayer}
    	    \end{scope}
         \begin{scope}[xshift=375pt]
	             \node [vertex] (a) at (-5/4,0) {};
    		    \node [vertex] (b1) at (0,-5/4) {};
	    	    \node [vertex] (b2) at (0,0) {};
		        \node [vertex] (c) at (5/4,0) {};
                \node [vertex] (d) at (0,5/4) {};
                \node [vertex] (e) at (0,10/4) {};
                \node[vertex] (f) at (0,-10/4) {};
        	
	        	\node [xshift=-\shft] at (b2) {$x$};
                     \node [xshift=\shft] at (b1) {$b$};
                    \node [xshift=-\shft] at (d) {$d$};
                    \node [xshift=-\shft] at (e) {$f$};
                    \node [yshift=\shft] at (a) {$a$};
                    \node [yshift=-\shft] at (c) {$e$};
                    \node [xshift=\shft] at (f) {$g$};
	        	
        		\begin{pgfonlayer}{bg}
	    	        \draw [edge] (a) -- (b1) -- (b2) -- (d) -- (c);
    	    	\draw [edge] (f) -- (b1);
	    	\draw [edge] (e) -- (d);        \draw [edge] (e) -- (d) -- (b2); \end{pgfonlayer}
    	    \end{scope}
	    \end{tikzpicture}
        \caption{The Hasse diagrams of the $\bigwedge$, $Y^d$, $S$, $S^+$, and $S^{++}$ posets.}
        \label{fig:WMposet}
    \end{center}
\end{figure}

For the proofs, we need to define five specific posets. See their Hasse diagrams in Figure~\ref{fig:WMposet}. Let
\begin{itemize}
    \item 
    $\bigwedge$ be the three element poset on $\{a,x,c\}$ with $a,x<c$, and $a,x$ unrelated;
    \item 
    $Y$ be the four element poset on $\{a,x,c,d\}$ with $a,x>c>d$ being all its cover relations;
    \item 
    $\cS$ be the five element poset on $\{a,b,x,d,e\}$ with $a>b<x<d>e$ being all its cover relations;
    \item 
    $\cS^+$ be the six element poset on $\{a,b,x,d,e,f\}$ with $a>b<x<d<f$ and $d>e$ being all its cover relations;
    \item 
    $\cS^{++}$ be the seven element poset on $\{a,b,x,d,e,f,g\}$ with $a>b$, $g<b<x<d<f$, and $d>e$ being all its cover relations.
\end{itemize}
    
\begin{proof}[Proof of Theorem \ref{thmtrue}]
Let $\eps>0$; note that it suffices to prove the theorem under the assumption that $\eps$ is sufficiently small and $n \in \mathbb N$ is sufficiently large.
Observe that if $x\in P$ is within distance at most 2 of all other elements of $P$, then a connected component of $P\setminus \{x\}$ is either an isolated vertex or contains a vertex $u$ that is contained in all edges of the Hasse diagram of that component, and where $u$ is adjacent to $x$ in the Hasse diagram of $P$. This and Remark \ref{equiv} imply that it is enough to prove the statement for the posets $\bigwedge,Y^d,\cS,\cS^+,\cS^{++}$ and $x$ as given in their definition and in Figure \ref{fig:WMposet}. For example, if $P$ is of height $2$ and $\gamma _1>0$, then there exists a $\gamma _2>0$ such that there is a choice of $y\in P$
so that $P(y,\gamma _2n) \subseteq \bigwedge (x,\gamma_1n)$ or  $P(y,\gamma _2n) \subseteq \bigwedge ^d (x,\gamma_1 n)$ for all sufficiently large $n\in \mathbb N$.

\medskip

 Consider first the poset $\bigwedge$; its height is 2.  It suffices to consider $\cF\subseteq 2^{[n]}$ of size exactly $(1+\varepsilon)\binom{n}{\lfloor n/2 \rfloor}$. 
 Let  $\delta:=\frac{\eps^2}{(1+\varepsilon)120}$.
As $n$ is large enough,  we can and will assume (here and in the cases of the four other posets as well) that all sets in $\cF$ have size between $n/2-n^{2/3}$ and $n/2+n^{2/3}$. 
Indeed, there are $o(\binom{n}{\lfloor n/2 \rfloor})$ sets in $2^{[n]}$ that do not have size  between $n/2-n^{2/3}$ and $n/2+n^{2/3}$, so we can ignore any such sets in $\cF$.
We say that $F\in \cF$ is of 
\begin{itemize}
    \item 
    \emph{Type 1} if there exist at least $\varepsilon n/30$ sets $G\in \cF$ with $|G|=|F|-1$ and $G\subset F$;
    \item 
    \emph{Type 2} if, for some $j \geq 2$, there exist at least $\varepsilon n^2/30$ sets $G\in \cF$ with $|G|= |F|-j$ and $G\subset F$;
    \item 
    \emph{Type 3} otherwise.
\end{itemize}

Let  $\cF^*$ denote  the subfamily of $\cF$ of sets of type 1 or 2.
Applying Lemma~\ref{fork} `upside down', we obtain that the subfamily $\cF\setminus \cF^*$ of sets of type 3 has size at most $(1+\varepsilon /2)\binom{n}{\lfloor n/2\rfloor}$. Therefore, $|\cF^*| \geq \varepsilon\binom{n}{\lfloor n/2\rfloor}/2$. Consider the  bipartite graph $B$ with parts $\cF$ and $\cF^*$, where there is an edge between $G\in \cF$ and $F\in \cF^*$ if and only if $G\subset F$. By definition of the types, the number of edges in $B$ is at least $\varepsilon n|\cF^*|/30 \ge \varepsilon^2 n\binom{n}{\lfloor n/2\rfloor}/60$; so there exists a set $F_x\in \cF$ with degree at least $2\delta n$ in $B$. $F_x$ will play the role of $x$ in the copy of $P(x,\delta n)$ that we will construct. At least $\delta n$ of $F_x$'s neighbours in $B$ are of the same type (i.e., all of type 1 or all of type 2); they will play the role of the many copies of $c$ in $P(x,\delta n)$. Let $\cF_c \subseteq \cF^*$ denote a family of such  sets so that 
$|\cF_c |=\delta n$. 

Suppose first that all sets in $\cF_c$ are of type 1. Each $F \in \cF_c$ has at least $\varepsilon n/30$ subsets in $\cF$ of size $|F|-1$; so there are at least $\eps n/30-|\cF_c|-1 \geq \eps n/60$ such subsets of $F$ that do not belong to $\cF_c\cup \{ F_x\}$. If $F,F'\in \cF_c$ have different sizes, then the subsets of $F$ and $F'$ that we consider are clearly distinct. If $F,F'$ have the same size, then they share at most one common subset of size $|F|-1$. Therefore, for all $F\in \cF_c$ we can pick $\eps {n}/{60}- |\cF_c|\geq \delta n$ distinct subsets of $F$ to obtain a copy of $
\bigwedge(x,\delta n)$ in $\mathcal F$.

Suppose next that $\cF_c$ consists only of sets of type 2. Then each of them contains at least ${\varepsilon n^2}/{30}-|\cF_c|-1 \geq \delta n |\cF_c|$ subsets that do not lie in $\cF_c \cup \{F_x\}$. So we can pick greedily $\delta n$ distinct subsets for each $F\in \cF_c$ to obtain a copy of $\bigwedge(x,\delta n)$. 

\medskip

Consider next the poset $Y^d$. 
 In this case it suffices to consider $\cF\subseteq 2^{[n]}$ of size exactly $(2+\varepsilon)\binom{n}{\lfloor n/2 \rfloor}$. Let $\delta_1$ be the output of the theorem on input $\eps/2$, $\bigwedge$ and $x$; thus, $\delta _1=\eps^2/ (480(1+\eps /2))$. 
 For $Y^d$ we show that we can take $\delta:= \delta_1/2$.
Now we say that $F\in \cF$ is of 
\begin{itemize}
    \item 
    \emph{Type 1} if there exist at least $\varepsilon n/30$ sets $G\in \cF$ with $|G|=|F|+1$ and $G\supset F$;
    \item 
    \emph{Type 2} if, for some $j\ge 2$, there exist at least $\varepsilon n^2/30$ sets $G\in \cF$ with $|G|= |F|+j$ and $G\supset F$;
    \item 
    \emph{Type 3} otherwise.
\end{itemize}
Applying Lemma~\ref{fork}, we obtain the subfamily of $\cF$ of sets of type 3 has size at most $(1+\varepsilon/2)\binom{n}{\lfloor n/2\rfloor}$. Therefore, the subfamily $\cF^* \subseteq \cF$ of sets of type 1 or 2 has size at least $(1+\varepsilon /2)\binom{n}{\lfloor n/2\rfloor}$. By the previous case, we obtain a copy of $\bigwedge(x,\delta _1 n)$ in $\cF ^*$. Let $F_x$ be the set corresponding to $x$, $\cF_c$ be the family of sets corresponding to the copies of $c$, and $\cF_a$ be the family of sets corresponding to copies of $a$. Suppose that there exist $F_c\in \cF_c, F_a\in \cF_a$ with $F_a\supset F_c$. Then we can exchange the roles of $F_a$ and $F_c$ and we still have a copy of $\bigwedge(x,\delta n)$. As after such a swap, the number of pairs $F_a\supset F_c$ decreases, after a finite number of changes, we can assume that no $F_a\in \cF_a$ contains any $F_c\in \cF_c$. 

Now we are ready to use our copy of $\bigwedge(x,\delta _1 n)$ to obtain a copy of $Y^d(x,\delta_1 n/2)$ by adding $\delta _1 n/2$ distinct supersets of each $F_c\in \cF_c$. We distinguish two cases according to the type of sets in $\cF_c$; at the price of working with $\bigwedge(x,\delta _1 n/2)$ rather than $\bigwedge(x,\delta _1 n)$, we can assume that all sets in $\cF_c$  are of the same type.

Suppose first that all sets in $\cF_c$ are of type 1. Then, by our assumption that no $F_a\in \cF_a$  contains any $F_c\in \cF_c$, any $F_c\in \cF_c$ has at least $({\varepsilon}/{30}-\delta_1/2)n$ supersets of size $|F_c|+1$ in $\cF$ that are not used in our copy of $\bigwedge(x,\delta _1 n/2)$. If $F_c$ and $F'_c$ have different sizes, then these supersets are distinct because they are also of different size. If $F_c,F'_c$ are of the same size, then they can share at most one common superset of size $|F_c|+1$; thus, as $\delta _1$ is much smaller than $\eps$, we can pick $\delta _1 n/2$ distinct supersets of each $F_c\in \cF_c$ to obtain $Y^d(x,\delta_1 n/2)$.
Suppose next that all $F_c\in\cF_c$ are of type 2. Then as ${\eps n^2}/{30}>\delta _1n|\cF_c|$, we can pick the supersets of each $F_c$ greedily to obtain a copy of $Y^d(x,\delta_1 n/2)$.

\medskip

Finally, we can consider the posets $S$, $S^+$, and $S^{++}$. As the proofs are almost identical, we 
only show the statement for $S^+$. Let  $\cF\subseteq 2^{[n]}$ be a family of size exactly 
$(3+\varepsilon)\binom{n}{\lfloor n/2\rfloor}$. 
Let $\delta_1$ be the  value of $\delta$ that we 
showed to exist for $\varepsilon/3$, $\bigwedge$ and $x$, and let $\delta_2$ be the value of $\delta$ 
that we showed to exist for $\varepsilon/3$, $Y^d$ and $x$. 
Set
$\delta^*:=\min \{\delta_1,\delta_2\}$.
Let $\cF_1$ be the subfamily of $\cF$ 
consisting of those sets $F\in \cF$ for which there exists no copy of $\bigvee(x,\delta^* n)\subset 
\cF$ in which $F$ plays the role of $x$. (The poset $\bigvee$ is the dual of $\bigwedge$.) Let $\cF_2$ 
be the subfamily of $\cF$ consisting of those sets $F\in \cF$ for which there exists no copy of 
$Y^d(x,\delta ^* n)\subset \cF$ in which $F$ plays the role of $x$. By definition of 
$\delta_1,\delta_2$ and $\delta ^*$, we have $|\cF_1|\le (1+\varepsilon/3)\binom{n}{\lfloor n/2\rfloor}$ and 
$|\cF_2|\le (2+\varepsilon/3)\binom{n}{\lfloor n/2\rfloor}$; so there exists a set $F_x\in \cF\setminus 
(\cF_1\cup \cF_2)$. By definition of $\cF_1$ and $\cF_2$, $F_x$ plays the role of $x$ both in a copy $\cG_1$ 
of $\bigvee(x,\delta^* n)$ and in a copy $\cG_2$ of $Y^d(x,\delta ^* n)$.  In $\cG_1\cup \cG_2$, the containment relations are as required as in a 
copy of $S^+(x,\delta^*n)$ with $F_x$ playing the role of $x$; the only problem that can occur is that 
$\cG_1$ and $\cG_2$ might overlap. We overcome this problem by decreasing the value of $\delta^*$ and 
filtering out sets from $\cG_1\cup \cG_2$.

For $y\in\{d,e,f\}$ let $\cG_y\subset \cG_2$ denote the 
family of sets corresponding to $y$ in $\cG_2$ and similarly for $y\in \{a,b\}$ let $\cG_y\subset 
\cG_1$ denote the family of sets corresponding to $y$ in $\cG_1$. Therefore, $\cG_d,\cG_e,\cG_f$ are 
pairwise disjoint and so are $\cG_a,\cG_b$. Also, $\cG_b\cap (\cG_d\cup \cG_f)=\emptyset$ as all sets 
in $\cG_b$ are proper subsets of $F_x$, while all sets in $\cG_d\cup \cG_f$ are proper supersets of 
$F_x$. For a set $G_b\in \cG_b$, we write $\cE(G_b)\subset\cG_a$ for the family of sets that correspond 
to the component of $\cG_1\setminus \{F_x\}$ containing $G_b$, and for a set $G_d\in \cG_d$ we write 
$\cE(G_d)\subset\cG_e\cup \cG_f$ for the family of sets that correspond to the component of 
$\cG_2\setminus \{F_x\}$ containing $G_d$.
\begin{itemize}
    \item 
    Consider a subfamily $\cG'_b\subset \cG_b$ with $|\cG'_b |=\frac{1}{2}|\cG_b|$. 
    Remove from $\cG_1$ (and thus $\cG_b$ and $\cG_a$ respectively) all elements of $\cG'_b$ as well as all elements from 
    $\{\cE(G'_b): G'_b\in \cG'_b\}$. 
    As we removed ${\delta^*n}/{2}$ sets from $\cG_b$, for every set $G_d\in \cG_d$ we have $|\cE(G_d)\cap \cG_b |\le {\delta^*n}/{2}$. We may therefore throw away exactly half of 
    $\cE(G_d) \cap \cG_e$ and exactly half of 
    $\cE(G_d) \cap \cG_f$ for every $G_d\in \cG_d$, so that now $\cE(G_d)$ is disjoint from $\cG_b$.
    We now must have that any overlaps between $\cG_1$ and $\cG_2$ occur between
     $\cG_a$ and $\cG_d\cup \cG_e\cup \cG_f$.
    Note that currently we have that $|\cG_b|=\delta ^*n/2$; $|\cE(G_b)|=\delta^* n$ for all
    $G_b \in \cG_b$; $|\cG_d |=\delta ^*n$; $|\cE(G_d) \cap \cG_e|=\delta^* n/2$ and 
    $|\cE(G_d) \cap \cG_f|=\delta^* n/2$ for all $G_d \in \cG_d$.
     
    \item 
    Consider a subfamily $\cG'_d\subset \cG_d$ with $|\cG'_d|=\frac{99}{100}|\cG_d|$. 
    Remove from $\cG_2$ (and thus $\cG_d$ and $\cG_e \cup \cG_f$ respectively) all elements of $\cG'_d$ as well as all elements from 
    $\{\cE(G'_d): G'_d\in \cG'_d\}$. 
    For every set $G_b\in \cG_b$ we therefore now have that $|\cE(G_b)\cap \cG_d |\le {\delta^*n}/{100}$.
    We can thus throw away exactly half of 
    $\cE(G_b) $ for every $G_b\in \cG_b$, so that now $\cE(G_b)$ is disjoint from $\cG_d$.
    We now must have that any overlaps between $\cG_1$ and $\cG_2$ occur between
     $\cG_a$ and $ \cG_e\cup \cG_f$.  Note that currently we have that $|\cG_b|=\delta ^*n/2$; $|\cE(G_b)|=\delta^* n/2$ for all
    $G_b \in \cG_b$; $|\cG_d |=\delta ^*n/100$; $|\cE(G_d) \cap \cG_e|=\delta^* n/2$ and 
    $|\cE(G_d) \cap \cG_f|=\delta^* n/2$ for all $G_d \in \cG_d$.
     
    \item We say that a $G_b \in \cG_b$ is \emph{destroyed} if $|\cE(G_b) \cap (\cG_e \cup \cG_f)|\geq 99|\cE(G_b)|/100=99 \delta ^* n/200$. Notice that at most $99|\cG_b|/100$ $G_b \in \cG_b$ are destroyed. 
    Indeed, suppose not. Then recalling that the sets $\cE(G_b)$ for $G_b \in \cG_b$ are disjoint, we obtain that
    $$
    |\cG_e \cup \cG_f|\geq \sum _ {G_b \in \cG_b} |\cE(G_b) \cap (\cG_e \cup \cG_f)|
    \geq \frac{99|\cG_b|}{100} \times \frac{99\delta ^* n}{200} =\left ( \frac{99\delta ^* n}{200}
    \right ) ^2.
    $$
    This is a contradiction though as
    $$  |\cG_e \cup \cG_f| =|\cG_d|\delta ^* n = \frac{(\delta ^*n)^2}{100}< \left ( \frac{99\delta ^* n}{200}
    \right ) ^2.$$
    Remove all destroyed $G_b$ from $\cG_b$ and also remove from $\cG_a$ everything from 
    $\{\cE(G_b) : G_b \text{ is destroyed}\}$.
    So now $|\cG_b|\geq \delta^* n/200$. For each $G_b \in \cG_b$ delete from  $\cE(G_b)$ all elements from 
    $\cE(G_b) \cap (\cG_e \cup \cG_f)$. Since each $G_b \in \cG_b$ is not destroyed, we still have that $|\cE(G_b)|=\delta^* n/200$. Recall also that 
    $|\cG_d |=\delta ^*n/100$, and
    $|\cE(G_d) \cap \cG_e|=\delta^* n/2$ and 
    $|\cE(G_d) \cap \cG_f|=\delta^* n/2$ for all $G_d \in \cG_d$. 
    
    We have ensured that $\cG_a$ and $ \cG_e\cup \cG_f$ are disjoint, and thus $\cG_1$ and $\cG_2$ are disjoint. Therefore, $\cG_1 \cup \cG_2$ contains a copy of $\cS^+(x,{\delta n})$ where $\delta:= \delta ^*/200$, as desired.
\end{itemize}


\end{proof}

\begin{proof}[Proof of Lemma \ref{lem:special}]
    The proof closely follows that of  Theorem~\ref{thmtrue}, but instead of using Lemma~\ref{fork} we apply Lemma \ref{fork+}. Again, as in the case of Theorem~\ref{thmtrue}, by Remark~\ref{equiv}, it is enough to only prove the result for $P=\bigwedge,Y^d,S,S^+,S^{++}$, but with the slightly stronger conclusion that we seek a copy of $P(x,n^{1.91})$. As  $x$ is at distance at most $2$ from any element of $P$,
    every element $y\in P$ has at most $n^{3.82}\ll  n^4$ sets corresponding to $y$ in a copy of $P(x,n^{1.91})$. 
    This combined with the fact that
we now consider families $\mathcal F$ that are four times as large as those in Theorem~\ref{thmtrue} means the proof is actually cleaner  than that of Theorem~\ref{thmtrue}. As such, we only explicitly show the case of $\bigwedge$.

    Let $0<\varepsilon<1/4$ be fixed, $n$ be sufficiently large  and $\cF\subseteq 2^{[n]}$ be of size 
    $(4+4\varepsilon)\binom{n}{\lfloor n/2\rfloor}=4(h(\bigwedge)-1+\varepsilon)\binom{n}{\lfloor n/2\rfloor}$. 
    Let $\cF^-$ be the family of sets $F\in \cF$ that do not contain ${\varepsilon}n^4/500$ other 
    sets of $\cF$ of size $|F|-j$ for any $j\ge 4$. By Lemma \ref{fork+}, 
    $|\cF^-|\le (4+\frac{4\varepsilon}{5})\binom{n}{\lfloor n/2\rfloor}$ and thus $\cF^+:=\cF\setminus \cF^-$ 
    has size at least $3{\varepsilon}\binom{n}{\lfloor n/2\rfloor}$. Consider the bipartite graph 
    $B$ with parts $\cF^+$ and $\cF$ such that $F^+\in \cF^+$ and $F\in \cF$ are joined by an edge if and 
    only if $F\subset F^+$. By definition of $\cF^+$, the number of edges in $B$ is at least 
    ${\varepsilon}n^4|\cF^+|/500\ge {\varepsilon^2}n^4\binom{n}{\lfloor n/2\rfloor}/200$; 
    therefore there exists $F_x\in \cF$ with degree at least ${\varepsilon^2}n^4/1000$ in $B$. $F_x$ will play 
    the role of $x$ in the copy of $\bigwedge(x,n^{1.91})$ that we will construct, and $n^{1.91}$ of its neighbours $F_c$ in $B$ will play the role 
    of $x$'s neighbours in $\bigwedge(x,n^{1.91})$. For each  of these  $n^{1.91}$ neighbours $F_c$ of $F_x$, 
    we  greedily pick $n^{1.91}$ distinct subsets of $F_c$ to obtain a copy of  
    $\bigwedge(x,n^{1.91})$ in $\mathcal F$.
\end{proof}


\subsection{Proof of Theorem~\ref{thm:randomgeneral}}

To prove Theorem~\ref{thm:randomgeneral} we need to apply the following container result.
\begin{lemma}\label{lem:con2}
  Let $P$ be any tree poset of height $h \leq 5$ and radius at most $2$.  Let $0<\eps<1/5$ and 
  let $\delta >0$ be as in Conjecture~\ref{conj1}.\footnote{Note that in Theorem~\ref{thmtrue} we proved that 
  Conjecture~\ref{conj1} holds for such $P$ and some $x \in P$; so we can indeed select $\delta$ as in  Conjecture~\ref{conj1}.}
  Suppose $n \in \mathbb N$ is sufficiently large and write $m:= \binom{n}{\lfloor n/2\rfloor }$.
Then there exist  functions 
$f : \binom{2^{[n]}}{\leq |P| n^{-1.9} 2^n} \to \binom{2^{[n]}}{\leq 4(h-1+\eps)m}$  and 
$g : \binom{2^{[n]}}{\leq |P|(4h-3)m/ \delta n} \to \binom{2^{[n]}}{\leq (h-1+\eps)m}$ 
such that 
for any $P$-free family $\mathcal F$ in $2^{[n]}$ there are disjoint subfamilies $\mathcal H_1, \mathcal H_2 \subseteq \mathcal F$ so that:
\begin{itemize}
    \item[(i)] $\mathcal H_1 \in \binom{2^{[n]}}{\leq |P| n^{-1.9} 2^n}$ and 
    $\mathcal H_1  \cup \mathcal H_2 \in \binom{ 2^{[n]}}{ \leq |P|(4h-3)m/ \delta n}$;
    \item[(ii)] $\mathcal H_1 \cup \mathcal H_2 $ and $g(\mathcal H_1 \cup \mathcal H_2)$ are disjoint;
    \item[(iii)] $\mathcal H_2 \subseteq f(\mathcal H_1)$;
    \item[(iv)] $\mathcal F \subseteq \mathcal H_1 \cup \mathcal H_2 \cup g(\mathcal H_1 \cup \mathcal H_2)$.
\end{itemize}

\end{lemma}
 Similarly to before we refer to the family $\mathcal H_1 \cup \mathcal H_2 \cup g(\mathcal H_1 \cup \mathcal H_2)$ produced by Lemma~\ref{lem:con2} as a \emph{container}; we call $\mathcal H_1$ and $\mathcal H_2$ the \emph{fingerprints}.

The proof of Lemma~\ref{lem:con2} closely follows that of 
Lemma~\ref{lem:con} except that we analyse the container algorithm over two separate stages. The idea of a multi-stage analysis of the graph container algorithm was first introduced in~\cite{bmt}. In particular, as in the proof of Theorem~\ref{randomsperner} presented in~\cite{bmt}, using a two-stage container algorithm will allow us to employ more careful calculations in the proof of Theorem~\ref{thm:randomgeneral}, thereby allowing us to deal with probabilities $p$ `close' to $1/n$.

\begin{proof}[Proof of Lemma~\ref{lem:con2}]
Let $P$ be as in the statement of the lemma and set $p:=|P|$. Let $x \in P$ be of distance at most two from every other element in the Hasse diagram of $P$.
Given any $0<\eps<1/5$, let $\delta >0$ be as in Conjecture~\ref{conj1}. Suppose $n \in \mathbb N$ is sufficiently large and let $m:= \binom{n}{\lfloor n/2\rfloor }$. 

Fix a total order $\mathcal O_{2^{[n]}}$ of the elements of ${2^{[n]}}$.
Let $\mathcal O_{P}$ be a total order of the vertices of $P$ such that the first vertex is $x$; the next vertices are the children of $x$; the next vertices are those vertices of distance two from $x$ in the undirected Hasse diagram of $P$, and so forth. 
We write $y_1,\dots, y_{p}$ for the vertices of $P$ ordered as in $\mathcal O_{P}$; so $y_1=x$.
Further, let $\mathcal P_{\text{blow}}$ denote the set of all copies of $P(x,t)$ in $2^{[n]}$ for all $t \in \mathbb N$. Let $\mathcal O_{\text{blow}}$ be a total order of the elements of $\mathcal P_{\text{blow}}$.

We now run our modified version of the graph container algorithm. The input of the algorithm is a $P$-free family $\mathcal F \subseteq 2^{[n]}$. The algorithm will output  fingerprints $\mathcal H_1$, $\mathcal H_2$  and a container $\mathcal H_1 \cup \mathcal H_2 \cup g(\mathcal H_1 \cup \mathcal H_2)$ where 
$\mathcal H_1 \cup \mathcal H_2 \subseteq \mathcal F \subseteq \mathcal H_1 \cup \mathcal H_2 \cup g(\mathcal H_1 \cup \mathcal H_2)$. 
We proceed in two stages.

\smallskip

{\bf \noindent Stage~$1$:}
Initially we set $\mathcal G^0_1:= 2^{[n]}$ and $\mathcal H^0_1:=\emptyset$. We will add vertices from $2^{[n]}$ to $\mathcal H^0_1 $ and remove vertices from $\mathcal G^0_1$ through the following iterative process, beginning at Step~$1$.

At Step~$i$, let $P_1(x,t) \in \mathcal P_{\text{blow}}$ be a copy of $P(x,t)$ in 
$\mathcal G^{i-1}_1\subseteq 2^{[n]}$ where we choose $t$ to be as large as possible. If there is more than one copy of $P(x,t)$ in $\mathcal G^{i-1}_1$ for this choice of $t$ we choose $P_1(x,t)$ to be the copy of $P(x,t)$ appearing earliest in the total order $\mathcal O_{\text{blow}}$. Let $x_1$ be the vertex of $P_1(x,t)$ that plays the role of $x$. 

\begin{itemize}
    \item { Suppose $x_1 \not \in \mathcal F$.} Then define 
    $\mathcal G^{i}_1:= \mathcal G^{i-1} _1\setminus \{x_1\}$ and $\mathcal H^i_1:=\mathcal H^{i-1}_1$. Proceed to Step~$i+1$.

    \item Suppose $x_1  \in \mathcal F$ and $t \geq  n^{1.9}$. We will update 
    $\mathcal G^{i-1}_1$ and $\mathcal H^{i-1}_1$ in several \emph{subphases} to obtain
    $\mathcal G^{i}_1$ and $\mathcal H^{i}_1$ respectively.
\begin{itemize}
 \item {\it Subphase~$1$.}   Delete $x_1$ from $\mathcal G^{i-1}_1$ and add it to $\mathcal H^{i-1}_1$. Proceed to Subphase~$2$.

 \item {\it Subphase~$2$.} 
 Recall that we write $y_1,\dots, y_{p}$ for the vertices of $P$ ordered as in $\mathcal O_{P}$. So $y_1=x$ and $y_2$ is a child of $x$ in $P$. In $P_1(x,t)$ there is a set of $t$ vertices $U_2 \subseteq \mathcal G^{i-1}_1$ so that each $x _2 \in U_2$ is adjacent to $x_1$ in the Hasse diagram of $P_1(x,t)$ and further the edge between
$x_1$ and $x_2$ is oriented in the same direction as the edge between $y_1(=x)$ and $y_2$ in the Hasse diagram for $P$. 

\smallskip

Look through the vertices  in $U_2\subseteq 2^{[n]}$ one by one, following the total order $\mathcal O_{2^{[n]}}$. If we run through the entire set $U_2$ without finding an element from $\mathcal F$ then we delete all $t$ vertices in $U_2$ from 
$\mathcal G^{i-1}_1$ and stop Step~$i$.

\smallskip

Otherwise, let $x_2 \in U_2$ denote the very first element from $\mathcal F$ that we discover in $U_2$. Add $x_2$ to $\mathcal H^{i-1}_1$ and remove from 
$\mathcal G^{i-1}_1$ both  $x_2$ and all elements of $U_2$ that  occur before $x_2$ in the total order
$\mathcal O_{2^{[n]}}$. So here we have deleted between $1$ and $t$ elements of
$\mathcal G^{i-1}_1$, and of these  only $x_2$ lies in $\mathcal F$. Proceed to Subphase~$3$.

\item {\it Subphase~$j$ (for $3\leq j \leq p$).} 
In the previous subphases we have defined $x_1,\dots, x_{j-1} \in P_1(x,t)$ corresponding to $y_1,\dots, y_{j-1} \in P$.

\smallskip

By definition of $\mathcal O_{P}$, the unique parent $y_k$ of $y_j$ in $P$ must be one of  $y_1,\dots, y_{j-1}$. Furthermore, by definition of $P(x,t)$,
there is a set of $t$ vertices $U_j\subseteq \mathcal G^{i-1}_1$ in $P_1(x,t)$ corresponding to $y_j$
so that  every  $x _j \in U_j$ is adjacent to $x_k$ in the Hasse diagram of $P_1(x,t)$, and moreover, the edge between
$x_k$ and $x_j$ is oriented in the same direction as the edge between $y_k$ and $y_j$ in the Hasse diagram for $P$. 

\smallskip

Look through the vertices  in $U_j\subseteq 2^{[n]}$ one by one, following the total order $\mathcal O_{2^{[n]}}$. If we run through the entire set $U_j$ without finding an element from $\mathcal F$ then we delete all $t$ vertices in $U_j$ from 
$\mathcal G^{i-1}_1$ and stop Step~$i$.

\smallskip

Otherwise, let $x_j \in U_j$ denote the very first element from $\mathcal F$ that we discover in $U_j$. Add $x_j$ to $\mathcal H^{i-1}_1$ and remove from 
$\mathcal G^{i-1}_1$ both $x_j$ and all elements of $U_j$ that  occur before $x_j$ in the total order
$\mathcal O_{2^{[n]}}$. So here we have deleted at most $t$ elements of
$\mathcal G^{i-1}_1$, of which only $x_j$ lies in $\mathcal F$. Proceed to Subphase~$j+1$.

\end{itemize}

\smallskip

At the end of these subphases relabel $\mathcal G^{i-1}_1$ as 
$\mathcal G^{i}_1$ and $\mathcal H^{i-1}_1$ as $\mathcal H^i_1$; then proceed to Step~$i+1$. 
Since $\mathcal F$ is $P$-free, for some $j \in [p]$, Subphase~$j$ must consider a set $U_j$ where $U_j \cap \mathcal F =\emptyset$. Thus, 
 in Step~$i$ we have deleted at least $t$ vertices from $\mathcal G^{i-1}_1$ to obtain $\mathcal G^i _1$ and have added at most $p-1$ vertices to 
$\mathcal H^{i-1}_1$ (at most one in each of the first $p-1$ subphases) to obtain $\mathcal H^i_1$.

    \item Suppose $x_1  \in \mathcal F$ and $t <  n^{1.9}$. Then define
     $\mathcal H_1:=\mathcal H^{i-1} _1 \cup \{x_1\}$ and $f(\mathcal H_1):= \mathcal G^{i-1}_1 \setminus \{x_1\}$ and proceed to Step~$1$ of Stage~$2$.
\end{itemize}

Note that in every step of Stage~$1$ we only add elements to $\mathcal H^{i-1}_1$ that lie in $\mathcal F$. So certainly $\mathcal H_1 \subseteq \mathcal F$. 
Similarly, by construction, $\mathcal F \subseteq \mathcal H_1 \cup f(\mathcal H_1)$, and $\mathcal H_1$ and $f(\mathcal H_1)$ are disjoint.
In every step of Stage~$1$ except the final step, if we add (at most $p-1$) elements to $\mathcal H^{i-1}_1$ we delete at least $ n^{1.9}$ elements from 
$\mathcal G^{i-1} _1 \subseteq 2^{[n]}$.  Therefore,
\begin{align}\label{eq1}
| \mathcal H_1|\leq \frac{(p-1)2^n}{  n^{1.9}}+1 \leq \frac{p \cdot 2^n}{  n^{1.9}}. 
\end{align}
Furthermore, by construction of Stage~$1$, the family $f(\mathcal H_1)\subseteq 2^{[n]}$ does not contain a copy of the $ n^{1.9}$-blow-up $P(x, n^{1.9})$. Thus, by Lemma~\ref{lem:special}, we have that 
$$|f(\mathcal H_1)|\leq 4(h-1+\eps ) m.$$

\smallskip

{\bf \noindent Stage~$2$.}
At the start of Stage~$2$ we set $\mathcal G^0_2:= f(\mathcal H_1)$ and $\mathcal H^0_2:=\emptyset$. We will add vertices from $f(\mathcal H_1)$ to $\mathcal H^0_2 $ and remove vertices from $\mathcal G^0_2$ through the following iterative process, beginning at Step~$1$.

At Step~$i$ of Stage~$2$, let $P_2(x,t) \in \mathcal P_{\text{blow}}$ be a copy of $P(x,t)$ in 
$\mathcal G^{i-1}_2\subseteq 2^{[n]}$ where we choose $t$ to be as large as possible. If there is more than one copy of $P(x,t)$ in $\mathcal G^{i-1}_2$ for this choice of $t$ we choose $P_2(x,t)$ to be the copy of $P(x,t)$ appearing earliest in the total order $\mathcal O_{\text{blow}}$. Let $x_1$ be the vertex of $P_2(x,t)$ that plays the role of $x$. 

\begin{itemize}
    \item { Suppose $x_1 \not \in \mathcal F$.} Then define 
    $\mathcal G^{i}_2:= \mathcal G^{i-1} _2\setminus \{x_1\}$ and $\mathcal H^i_2:=\mathcal H^{i-1}_2$. Proceed to Step~$i+1$.

    \item Suppose $x_1  \in \mathcal F$ and $t \geq  \delta n$. We then update 
    $\mathcal G^{i-1}_2$ and $\mathcal H^{i-1}_2$ in at most $p$ {subphases} to obtain
    $\mathcal G^{i}_2$ and $\mathcal H^{i}_2$ respectively.
    These subphases are identical to those subphases described in Stage~$1$, just with $\mathcal G^{i-1}_2$, $\mathcal H^{i-1}_2$ and $P_2(x,t)$ now playing the roles of 
    $\mathcal G^{i-1}_1$, $\mathcal H^{i-1}_1$ and $P_1(x,t)$
    respectively.

    At the end of the subphases relabel $\mathcal G^{i-1}_2$ as $\mathcal G^{i}_2$ and $\mathcal H^{i-1}_2$ as $\mathcal H^i_2$; then proceed to Step~$i+1$ of Stage~$2$. 
    Since $\mathcal F$ is $P$-free, for some $j \in [p]$, Subphase~$j$ must consider a set $U_j$ where $U_j \cap \mathcal F =\emptyset$. Thus, in Step~$i$ of Stage~$2$, we have deleted at least $t$ vertices from $\mathcal G^{i-1}_2$ to obtain $\mathcal G^i _2$ and have added at most $p-1$ vertices to $\mathcal H^{i-1}_2$ (at most one in each of the first $p-1$ subphases) to obtain $\mathcal H^i_2$.

 \item Suppose $x_1  \in \mathcal F$ and $t <  \delta n$. Then define
     $\mathcal H_2:=\mathcal H^{i-1} _2 \cup \{x_1\}$ and $g(\mathcal H_1 \cup \mathcal H_2):= \mathcal G^{i-1}_2 \setminus \{x_1\}$ and terminate the algorithm.
\end{itemize}

Note that in every step of Stage~$2$ we only add elements to $\mathcal H^{i-1}_2$ that lie in $\mathcal F \cap f(\mathcal H_1)$. So certainly $\mathcal H_2 \subseteq \mathcal F$. 
Further, as $\mathcal H_1 $ and $f(\mathcal H_1)$ are disjoint,
this implies $\mathcal H_1$ and $\mathcal H_2$ are disjoint.
Similarly, by construction,
conditions (ii)--(iv) of the lemma hold.

In every step of Stage~$2$ except the final step, if we add (at most $p-1$) elements to $\mathcal H^{i-1}_2$ we delete at least $\delta n$ elements from 
$\mathcal G^{i-1} _2 \subseteq 2^{[n]}$.  
Moreover, 
recall that $\mathcal G^0_2:= f(\mathcal H_1)$ and
$|f(\mathcal H_1)|\leq 4(h-1+\eps ) m.$
Therefore,
$$| \mathcal H_2|\leq \frac{(p-1)4(h-1+\eps ) m}{\delta  n}+1 \leq \frac{p \cdot 4(h-1+\eps ) m}{\delta  n}. $$
Combining this with (\ref{eq1}) we see that condition (i) of the lemma holds.

Note that $g(\mathcal H_1 \cup \mathcal H_2)\subseteq 2^{[n]}$ is defined so that it does not contain a copy of $P(x,\delta n)$. Thus, by Theorem~\ref{thmtrue},
$$|g(\mathcal H_1 \cup \mathcal H_2)| \leq (h-1+\eps)m. $$

Therefore all that remains to check is that the functions $f$ and $g$ are well-defined. That is, if for two $P$-free families 
$\mathcal F$ and $\mathcal F'$ the algorithm outputs the same $\mathcal H_1$, then $f(\mathcal H_1)$ is defined the same, and furthermore, if $\mathcal H_1 \cup \mathcal H_2$ is the same for both runs of the algorithm then so is 
$g(\mathcal H_1 \cup \mathcal H_2)$.
 This follows though immediately from the definition of the algorithm. 
Indeed, suppose one is only presented with the fingerprint $\mathcal H_1$ that is outputted by Step~$1$ of the algorithm. Then one can completely identify every action taken during Step~$1$ of the  algorithm. Similarly, if one is further given 
the second fingerprint $\mathcal H_2$, then one can 
completely identify every action taken during Step~$2$ of the algorithm.
\end{proof}
With Lemma~\ref{lem:con2} at hand, we can now prove Theorem~\ref{thm:randomgeneral}.

\begin{proof}[Proof of Theorem~\ref{thm:randomgeneral}.]
Let $P$ be any tree poset of height $h \leq 5$ and radius at most $2$.
Fix $\eps > 0$; it suffices to prove the theorem under the assumption that $\eps < 1/5$. Let $\eps _1 := \eps/4$ and let $\delta >0$ be as in Conjecture~\ref{conj1} on input $\eps _1$.
Define $C:=10^{10} |P|(4h-3) h^4 \eps ^{-5} \delta ^{-1}$. Let $p > C/n$ and $m:=\binom{n}{\lfloor n/2 \rfloor }$.

Note that the middle $h-1$ layers of $2^{[n]}$ form a $P$-free subfamily of $2^{[n]}$ 
of size $(h-1-o(1)) m$. 
So w.h.p. $\P (n,p)$ contains a $P$-free family of size at least
$(h-1-\eps ) p m$.
It remains to show that, w.h.p., $\P(n,p)$ contains no $P$-free family of size greater than $(h-1+\eps)pm$. We now follow the proof of Theorem~\ref{randomsperner} given in~\cite{bmt} closely.

Apply Lemma~\ref{lem:con2} with $\eps _1$ playing the role of $\eps$.
Suppose for a contradiction that $\P(n,p)$ does contain some $P$-free family $\mathcal F$ with $|\mathcal F| > (h-1+\eps)pm$. Let $\mathcal H_1$ and 
 $\mathcal H_2$ denote the fingerprints for $\mathcal F$ given by Lemma~\ref{lem:con2}. 
 Then since $\mathcal H_1 \cup \mathcal H_2 \subseteq \mathcal F$, we must have that $\mathcal H_1 \cup \mathcal H_2 \subseteq \P(n,p)$. Further, at least 
 $|\mathcal F| - |\mathcal H_1 \cup \mathcal H_2| \geq (h-1+\eps)pm -  |P|(4h-3)m/ \delta n\geq (h-1+\eps/2)pm$  elements of $g(\mathcal H_1 \cup \mathcal H_2)$ must be in 
 $\P(n,p)$.

Note that the number of possibilities for $\mathcal H_1$ is 
$\binom{2^n}{\leq |P|n^{-1.9} 2^n}$, and for each possibility 
the probability that $\mathcal H_1 \subseteq \P(n,p)$ is 
$p^{|\mathcal H_1|}$. For any fixed $\mathcal H_1$ we have 
$|f(\mathcal H_1)| \leq 4(h-1+\eps _1)m\leq (4h-3)m$ and 
$\mathcal H_2 \subseteq 
f(\mathcal H_1)$, so the number of possibilities for 
$\mathcal H_2$ is at most $\binom{(4h-3)m}{\leq |P|(4h-3)m / \delta n}$, and for each possibility the probability that $\mathcal H_2 \subseteq \P(n,p)$ is $p^{|\mathcal H_2|}$.

Furthermore, for any fixed $\mathcal H_1$ and $\mathcal H_2$ we have 
$g(\mathcal H_1 \cup \mathcal H_2) \leq (h-1+\eps_1) m = (h-1+\eps/4)m$, so the expected number of elements of $g(\mathcal H_1 \cup \mathcal H_2)$ 
selected for $\P(n,p)$ is at most $(h-1+\eps/4)pm$. By the Chernoff bound for the binomial distribution, the probability that at least $(h-1+\eps/2)pm$ elements of $g(\mathcal H_1 \cup \mathcal H_2)$ are selected for $\P(n,p)$ is therefore at most $e^{-\eps^2pm/(100h^2)}$. 

Taking a union bound, we conclude that the probability that $\P(n,p)$ contains a $P$-free family of size greater than $(h-1+\eps) pm$ is at most  
\begin{align*}
\Pi &:= \sum _{0\leq a \leq |P| n^{-1.9} 2^n}  \ 
\sum _{0\leq b \leq |P|(4h-3)m/\delta n } \binom{2^n}{a}\cdot p^{a}\cdot \binom{(4h-3)m}{b}\cdot p^{b} \cdot  
e^{-\eps^2 pm/(100h^2)} \\
&
\leq \left (|P|n^{-1.9} 2^n +1 \right) \left (|P|(4h-3)m/\delta n  +1 \right ) 
\binom{2^n}{|P|n^{-1.9} 2^n} \cdot p^{|P|n^{-1.9} 2^n}
\binom{(4h-3)m}{|P|(4h-3)m/\delta n }\\ & \ \ 
\cdot  p^{|P|(4h-3)m/\delta n } \cdot 
e^{-\eps^2 pm/(100h^2)}.
\end{align*}
Note that for large $n$, with room to spare we have
$$(|P|n^{-1.9} 2^n +1) (|P|(4h-3)m/\delta n  +1) \leq e^{\eps^2 pm/(400h^2)}$$ and
$$\binom{2^n}{|P|n^{-1.9} 2^n} \cdot p^{|P|n^{-1.9} 2^n} \leq e^{\eps^2 pm/(400h^2)}.$$
Further, as $C=10^{10} |P|(4h-3) h^4 \eps ^{-5} \delta ^{-1}$, for sufficiently large $n$ we have that
$$\binom{(4h-3)m}{|P|(4h-3)m/\delta n }
\cdot  p^{|P|(4h-3)m/\delta n } \leq e^{\eps^2 pm/(400h^2)}.$$
Therefore, the upper bound $\Pi$ on the probability is $o(1)$, as required.
\end{proof}

\section{Further discussion and observations}\label{sec:further}
For tree posets $P$ of radius at most $2$,  we have   asymptotically
determined the number of $P$-free families in $2^{[n]}$ (Corollary~\ref{cor1}) and have also resolved the random version of the $P$-free problem (Theorem~\ref{thm:randomgeneral}).
In general though, these questions for other posets $P$ remain wide open.
However, as we will now see, if one has solved either of these problems for some poset $P$, then one can often deduce further results from this.

Let $P_{t,h}$ denote the (upward) monotone $t$-ary tree poset of height $h$; the \emph{$t$-fork} $\vee_t$ is $P_{t,2}$.
Let $\mathcal D_k$ denote the \emph{$k$-diamond}; that is, the poset with a unique maximal element above all other elements, a unique minimal element below all other elements, and $k$ (incomparable) elements in between. Note that $\Diamond:=\mathcal D_2$  is the \emph{diamond}.

\begin{observation}\label{monotone}\

\begin{enumerate}
    \item 
    Given any $t, h \in \mathbb N,$ and $\eps>0$, there exists $C>0$ such that the following holds. If $p>C/n$ then w.h.p. the largest $P_{t,h}$-free family in $\mathcal P(n,p) $ has  size
$(h-1 \pm \eps )p\binom{n}{\lfloor n/2 \rfloor}.$ 
    \item 
    The number of $P_{t,h}$-free families in $2^{[n]}$ is $2^{(h-1+o(1))\binom{n}{\lfloor n/2\rfloor}}$.
\end{enumerate} 
\end{observation}

\begin{proof}
    The lower bound of (1) follows by noting that, w.h.p.,  $\mathcal P(n,p) $ contains at least 
    $(h-1 - \eps )p\binom{n}{\lfloor n/2 \rfloor}$ elements from the $h-1$ middle layers of $2^{[n]}$.
    The lower bound of (2) can be seen by considering all subfamilies of the $h-1$ middle layers of $2^{[n]}$.
    
    The proof of the upper bound of (1) proceeds by induction on $h$, with the base case $h=2$ proved by Hogenson \cite{hog}. Suppose (1) holds for any choice of the height of the tree less than $h$ where $h \geq 3$. Fix any $t\in \mathbb N$ and $\eps >0$.
    By the induction hypothesis, there exists $C_1>0$ such that if $p>C_1/n$ then, w.h.p, the largest 
    $P_{2t^h,2}$-free family in $\mathcal P(n,p) $ has size at most $ (1+\varepsilon/2)p\binom{n}{\lfloor n/2\rfloor}$. Further, 
    there exists $C_2>0$ such that if $p>C_2/n$ then, w.h.p, the largest 
    $P_{t,h-1}$-free family in $\mathcal P(n,p) $ has size at most $ (h-2+\varepsilon/2)p\binom{n}{\lfloor n/2\rfloor}$.

    Let $C:=\max\{C_1,C_2\}$ and consider $p>C/n$.
    Consider any $P_{t,h}$-free family $\cF\subseteq \P(n,p)$. Define $\cF':=\{F\in \cF: |\cU (F)\setminus \{F\}|<2t^{h}\}$. As $\cF'$ is  $P_{2t^h,2}$-free then, w.h.p., 
    $|\cF'|\le (1+\varepsilon/2)p\binom{n}{\lfloor n/2\rfloor}$.
    
    Furthermore,  observe that $\cF\setminus \cF'$ is $P_{t,h-1}$-free, as otherwise a copy of $P_{t,h-1}$ in $\cF\setminus \cF'$ could be extended greedily to a copy of $P_{t,h}$ in $\cF$. Indeed, any set $F$ corresponding to a leaf of $P_{t,h-1}$ is contained in at least $2t^h$ many other sets of $\cF$, at least $t^h$ of which do not belong to the copy of $P_{t,h-1}$. As there are $t^{h-2}$ leaves, we can pick pairwise disjoint families of $t$ supersets for every leaf. 
    
    Thus, w.h.p.,  $ |\cF\setminus \cF'|\leq (h-2+\varepsilon/2)p\binom{n}{\lfloor n/2\rfloor}$, yielding $|\cF|\le (h-1+\varepsilon)p\binom{n}{\lfloor n/2\rfloor}$ w.h.p., as desired.

    The upper bound of (2) follows similarly by induction on $h$ with the base case $h=2$ proved in \cite{hog}. Any $P_{t,h}$-free family $\cF\subseteq 2^{[n]}$ can be partitioned into $\cF'$ and $\cF\setminus \cF'$ as above; so the number of $P_{t,h}$-free families in $2^{[n]}$
    is at most the product of the number of $P_{t,h-1}$-free families and the number of $P_{2t^h,2}$-free families. The statement then follows by induction.
\end{proof}
Since any upward monotone tree poset $P$ of height $h$ is contained in $P_{t,h}$ for $t$ large enough, Observation~\ref{monotone} holds for any such $P$. Moreover, it holds for a somewhat wider class of tree posets. We say that $P$ is a \emph{two-way monotone tree poset} if there exists $x\in P$ such that $\cU(x)\cup\cD(x)=P$ where  $\cU(x)$ is the set of elements of $P$ that are supersets of $x$ and $\cD(x)$ is the set of elements of $P$ that are subsets of $x$. It can be shown that Observation~\ref{monotone} holds for two-way monotone tree posets of height $h$ as well.

\smallskip

In this paper we have considered the problem of determining for what probabilities $p$  the largest $P$-free family in
$\P(n,p)$ has size $p \cdot(1+ o(1)) La(n,P)$ with high probability. 
It is of course natural to consider other values of $p$ (i.e., when the 
largest $P$-free family in
$\P(n,p)$ has size more than $p \cdot (1+ o(1))  La(n,P)$); see~\cite{bmt, osthus}.
    The following is such a result for the diamond.
    \begin{observation}
  Given $\eps>0$ there exist $C_1,C_2>0$ such that the following holds. If $C_1/n< p< C_2 /n^{2/3}$ then w.h.p. the largest $\Diamond$-free family in $\mathcal P(n,p) $ has size  
$(3 \pm \eps )p\binom{n}{\lfloor n/2 \rfloor}.$ 
\end{observation}

\begin{proof}
    The upper bound follows from the fact that a $\Diamond$-free family $\cF \subseteq \P(n,p)$ cannot contain a 4-chain which is a special copy of $\Diamond$. Therefore $\cF$ can be partitioned into three antichains $\cF_1,\cF_2,\cF_3$; thus the upper bound on $|\mathcal F|$ follows from Theorem~\ref{randomsperner} applied to all $\cF_i$s.

    To see the lower bound, observe that the expected number $E_1$ of sets in the middle three layers of $\P(n,p)$ is $(3\pm \varepsilon/4)p\binom{n}{\lfloor n/2\rfloor}$. Also, the expected number $E_2$ of copies of $\Diamond$ in the middle three layers of $\P(n,p)$ 
    is $(1\pm \varepsilon/4)p^4\frac{n^2}{8}\binom{n}{\lfloor n/2\rfloor}$. Using the Markov and Chebyshev inequalities, it follows that if $C_2>0$ is small enough, then w.h.p. $E_1\ge (3-\varepsilon/2)p \binom{n}{\lfloor n/2\rfloor}$ and $E_2\le (1+ \varepsilon/2)p^4\frac{n^2}{8}\binom{n}{\lfloor n/2\rfloor}\le \frac{\varepsilon}{2}p\binom{n}{\lfloor n/2\rfloor}$. So removing one element from each copy of a $\Diamond$   leaves us with a $\Diamond$-free family of size at least $(3-\varepsilon)p \binom{n}{\lfloor n/2\rfloor}$ in $\P(n,p)$.
\end{proof}

As Theorem~\ref{thm:count} articulates, Conjecture~\ref{conj1} has implications to the counting problem for tree posets $P$. The resolution of Conjecture~\ref{conj1} would also likely allow one to generalise Theorem~\ref{thm:randomgeneral} to all tree posets $P$ of height $h$, not just those of radius $2$. For this one would also need an analogue of Lemma~\ref{lem:special}.

One might wonder why our proof of Theorem~\ref{thmtrue} does not work for posets of radius larger than 2, and why Bukh's original argument~\cite{bukh} cannot be applied to prove Theorem~\ref{thmtrue} for arbitrary tree posets and arbitrary element $x$. Let us consider first the latter question. Bukh's proof uses some preliminary structural results on tree posets. Namely, he shows that any tree poset of height $h$ is a subposet of a height $h$ \textit{saturated} tree poset, i.e., one in which all maximal chains contain $h$ elements. Then he shows that any saturated tree poset $T$ of height $h$ can be obtained from $C_h$ by `adding intervals' $I_1,I_2,\dots,I_m$, i.e., branches of the Hasse diagram. This enables him to obtain an inductive proof by showing that $T\setminus I_m$ can be embedded not only into any large enough family $\cF$, but to nice, extendable sets of $\cF$, where nice means that a \textit{constant} number of sets cannot ruin the extending property. There are two points of this proof that we do not see how to save to generalise Theorem~\ref{thmtrue}. As a height $h$ tree poset can have arbitrary radius, it is not enough to extend the embedding of the blow-up by one interval, but by a polynomial (of arbitrary high degree!) number of   blow-ups of intervals, and now the extension property should not be ruined by  a polynomial number of sets (rather than a constant number). Even more importantly, we do not know how to start the induction. Bukh's result for chains (even with the extension property) is immediate from any supersaturation result for $h$-chains, but we are unable to prove Theorem~\ref{thmtrue} for $C_4$ and $x$ being its minimal element. 

This brings us to the other question: why our proof does not work for the blow-up of trees with larger radius. Getting the result for $\bigvee$ and $\bigwedge$ is basically an exercise if one is familiar with the notion of Lubell-mass. Then for the poset $Y$, we try to apply some kind of an induction: we find lots of blow-ups of $\bigvee$ and then a blow-up of $\bigvee$ in the root of the blow-ups of the first round. All we need to take care of are the overlaps of these blow-ups. This is done via introducing type 1, 2, and 3 sets and observing that if we work with type 1 sets, then the pairwise overlaps are of size at most 1 and so at the price of shrinking $\delta$, we get our pairwise disjoint blow-ups, while if we work with type 2 sets, then we need only a linear number of sets, while the blow-ups are quadratic. The latter could work for arbitrary long chains, but type 1 sets are problematic: type 1 copies of blow-ups of $Y$ can intersect in a linear number of sets, and we were unable to handle how to get many pairwise disjoint copies of them.

Because of the above reasoning, the full resolution of each of Conjecture~\ref{conj1}, Conjecture~\ref{conj:count} and the random version of the $P$-free problem does currently seem out of reach.  It would therefore be interesting to resolve these problems for other natural classes of posets $P$. For example, it would be interesting to extend our results to cover all tree posets $P$ of height $2$, particularly those posets whose (undirected) Hasse diagram is a path. In~\cite{griggs3}, $La(n,\mathcal D_k)$ was asymptotically determined for infinitely many choices of $k$. It would be interesting to resolve  Conjecture~\ref{conj:count} and the random version of the $P$-free problem for such $\mathcal D_k$. It is also natural to seek induced versions of our results, including Corollary~\ref{cor1} and Theorem~\ref{thm:randomgeneral}.

\subsection*{Acknowledgments}
The research in this paper was carried out during the \emph{Focused Workshop on Saturation in the Forbidden Subposet Problem} held at the Erd\H{o}s Center, Budapest (July 31st to August 4th, 2023). The authors are grateful to the Erd\H{o}s Center for the nice working environment and  to Maria-Romina Ivan, Ryan Martin, D\'aniel Nagy and Casey Tompkins for helpful discussions during this workshop. We are also grateful to Ryan Martin for comments on this manuscript, and to  the referees for their helpful and careful reviews.

The first author was supported by NKFIH grant   FK 132060.
The second author was  supported by EPSRC  grant EP/V002279/1.
\smallskip

{\noindent \bf Data availability statement.}
There are no additional data beyond that contained within the main manuscript.


\end{document}